\documentclass[12pt]{article}

\usepackage[margin=1in]{geometry}
\usepackage{charter}
\usepackage[noblocks]{authblk}

\usepackage{amsmath,amsfonts,bm}

\newcommand{\vect}[1]{\mathbf{#1}} 

\def\1{\bm{1}}

\def\bxi{\boldsymbol{\xi}}
\def\bpi{\boldsymbol{\pi}}

\def\NN{\text{NN}}

\def\argmax{\mathop{\rm arg\,max}}%
\def\argmin{\mathop{\rm arg\,min}}%

\DeclareMathAlphabet{\mathsfit}{\encodingdefault}{\sfdefault}{m}{sl}
\SetMathAlphabet{\mathsfit}{bold}{\encodingdefault}{\sfdefault}{bx}{n}

\def\gX{{\mathcal{X}}}
\def\gY{{\mathcal{Y}}}

\newcommand{\R}{\mathbb{R}}

\usepackage[usenames,dvipsnames]{xcolor}
\usepackage{graphicx, pgfplots, standalone, pgf, tikz}
\usepackage[inline]{enumitem}
\usepackage[justification=justified,singlelinecheck=false]{caption}
\usepackage{rotating}
\usepackage{subfig}

\usepackage{multirow}
\usepackage{multicol}

\usepackage{algorithmic}
\newcommand\PCref[1]{(\Cref{#1})}
\usepackage{booktabs}

\usepackage[sort&compress,numbers]{natbib}

\usepackage{comment, url,float, setspace, xspace}

\usepackage{amsmath, amssymb, amsthm, amsfonts, mathtools, array, bm, xfrac}

\usepackage[ruled,vlined,linesnumbered]{algorithm2e}

\usepackage[short,c2,nocomma]{optidef}

\usepackage[colorlinks,
            linkcolor=BlueViolet,
            citecolor=blue,
            urlcolor=blue]{hyperref}

\usepackage[nameinlink]{cleveref}

\usepackage{csvsimple, booktabs, longtable, multirow, multicol, adjustbox}	

\usepackage[xindy, acronyms, nowarn]{glossaries}   
\makeglossaries

\definecolor{lime}{HTML}{A6CE39}
\DeclareRobustCommand{\orcidicon}{
	\begin{tikzpicture} \draw[lime, fill=lime] (0,0) circle [radius=0.16] node[white] { {\fontfamily{qag}\selectfont \tiny ID} };
	\draw[white, fill=white] (-0.0625,0.095) circle [radius=0.007];
	\end{tikzpicture} \hspace{-2mm}
}

\foreach \x in {A, ..., Z}{\expandafter\xdef\csname orcid\x\endcsname{\noexpand\href{https://orcid.org/\csname orcidauthor\x\endcsname}{\noexpand\orcidicon} } }
\newcommand{\converged}{\textit{converged}}
\newcommand{\true}{\textbf{True}}
\newcommand{\false}{\textbf{False}}

\newcommand{\method}{\textsc{Neur2RO}}
\newcommand{\methodmip}{$\textsc{N}_{\textsc{MILP}}$}
\newcommand{\methodipo}{$\textsc{N}_{\textsc{IPO}}$}
\newcommand{\methodngc}{$\textsc{N}_{\textsc{NGC}}$}
\newcommand{\static}{SRO}

\newtheorem{theorem}{Theorem} 

\crefname{line}{Line}{Lines}
\crefname{lemma}{Lemma}{Lemmata}
\crefname{theorem}{Theorem}{Theorems}
\crefname{proposition}{Proposition}{Propositions}
\crefname{algorithm}{Algorithm}{Algorithms}
\crefname{observation}{Observation}{Observations}
\crefname{equation}{}{}
\crefname{definition}{Definition}{Definitions}
\crefname{claim}{Claim}{Claim}
\crefname{corollary}{Corollary}{Corollaries}
\crefname{remark}{Remark}{Remarks}
\crefname{example}{Example}{Examples}
\crefname{figure}{Figure}{Figures}
\crefname{section}{Section}{Sections}
\crefname{table}{Table}{Tables}

\usepackage{hhline}

\title{Deep Learning for Two-Stage\\Robust Integer Optimization}

\author[1,2]{Justin Dumouchelle\footnote{Corresponding author: \href{mailto:justin.dumouchelle@mail.utoronto.ca}{justin.dumouchelle@mail.utoronto.ca}}}
\author[3]{Esther Julien}
\author[4]{Jannis Kurtz}
\author[1,2]{\\Elias B. Khalil}

\affil[1]{University of Toronto, Canada}
\affil[2]{SCALE AI Research Chair in Data-Driven Algorithms for Modern Supply Chains, Canada}
\affil[3]{TU Delft, The Netherlands}
\affil[4]{University of Amsterdam, The Netherlands}

\begin{document}
\maketitle
\begin{abstract}
\noindent Robust optimization is an established framework for modeling optimization problems with uncertain parameters. While static robust optimization is often criticized for being too conservative, two-stage (or adjustable) robust optimization (2RO) provides a less conservative alternative by allowing some decisions to be made after the uncertain parameters have been revealed. Unfortunately, in the case of integer decision variables, existing solution methods for 2RO typically fail to solve large-scale instances, limiting the applicability of this modeling paradigm to simple cases. We propose \method{}, a deep-learning-augmented instantiation of the column-and-constraint-generation (CCG) algorithm, which expands the applicability of the 2RO framework to large-scale instances with integer decisions in both stages. A custom-designed neural network is trained to estimate the optimal value and feasibility of the second-stage problem. 
The network can be incorporated into CCG, leading to more computationally tractable subproblems in each of its iterations. The resulting algorithm enjoys approximation guarantees which depend on the neural network's prediction error. In our experiments, \method{} produces high-quality solutions quickly, outperforming state-of-the-art methods on two-stage knapsack, capital budgeting, and facility location problems. Compared to existing methods, which often run for hours, \method{} finds better solutions in a few seconds or minutes. Our code is available at \url{https://github.com/khalil-research/Neur2RO}.
\end{abstract}

\section{Introduction} 
\label{sec:intro}
A wide range of real-world optimization problems in logistics, finance, and healthcare, among others, are modeled by discrete optimization models \citep{petropoulos2023operational}. While mixed-integer (linear) problems (MILP) can be challenging to solve, the size of the problems that modern MILP solvers can handle has increased significantly thanks to algorithmic improvements \citep{wolsey2020integer,achterberg2013mixed}. In recent years, machine learning (ML) has emerged as a tool for supporting classical MILP solvers in many different ways~\citep{zhang2023survey,bengio2021machine}. Based on similar instances from the same class of problems, ML models are trained to make search-related decisions such as branching variable selection, or predict (partial) solutions to the problem. While most of the work in this vein has focused on deterministic problems, decision-makers often deal with uncertainty in some problem parameters, e.g., due to forecasting or measurement errors in quantities of interest such as customer demands in transportation and production problems. Besides the stochastic optimization approach, for which learning-based heuristics have been proposed recently~\citep{dumouchelle2022neur2sp,alcantara2024quantile,wuhgcn2sp,kronqvist2023alternating,bae2023deep,chan2022machine,larsen2024fast}, another popular way to incorporate uncertainty into optimization models is \textit{robust optimization}, where the goal is to find solutions that are optimal for the worst-case realization of uncertain parameters that belong to a pre-defined uncertainty set \citep{ben2009robust}. For applications where some of the decisions can be made on the fly after the uncertain parameters are realized (e.g., slack variables or assignment decisions), the static robust optimization approach can lead to conservative solutions. A better model in this situation is two-stage robust optimization (2RO) \citep{ben2004adjustable, yanikouglu2019survey}, where it is assumed that some of the decisions have to be made right now (here-and-now decisions), while some decisions can be made later after the uncertain parameters are known (wait-and-see decisions). The goal is to find here-and-now decisions that are optimal in the worst-case over all scenarios in the given uncertainty set, where, for each scenario, the best wait-and-see decision is implemented.

In this work, we consider the setting in which the decision-maker routinely solves instances of the same 2RO problem that differ in some parameter values. We propose~\method, a deep-learning-augmented instantiation of the column-and-constraint generation (CCG) algorithm for 2RO.  CCG is an exact algorithm for 2RO problems that alternates between solving a main problem to produce a candidate here-and-now solution that is robust for a subset of scenarios from the uncertainty set, and solving a bilevel problem to generate a violating scenario from the uncertainty set that cuts off said solution.  CCG often requires a large amount of computation time due to its high per-iteration costs. The goal of this paper is to solve both subproblems by fast neural-network-supported \textit{surrogate} problems, leading to fast and good-quality solutions for 2RO.
Historical data---in the form of previously encountered instances---can be leveraged to train the {neural network}. A shorter version of this manuscript is published in~\citep{dumouchelle2024neurro}; this version extends~\citep{dumouchelle2024neurro} to include constraint uncertainty along with additional theoretical and experimental results.

\paragraph{Contributions.} 
\begin{enumerate}[label=(\roman*)]
    \item We present \method{}, one of the first ML approaches to 2RO. Our method uses a trained neural network to predict the optimal value and the feasibility of the second-stage problem of the 2RO problem.
    \item We show that the trained network can be incorporated into the classical CCG framework in a non-trivial fashion, leading to significantly more tractable subproblems that have to be solved during the CCG execution.
    \item We develop a set-based neural network architecture that can embed instances of 2RO of arbitrary size, which leads to high-quality predictions while having a relatively small number of parameters. Moreover, the design of our architecture is tailored to the structure of the CCG subproblems, resulting in tractable MILP representations of the networks.  
    \item We show that our algorithm has finite convergence and achieves an additive approximation guarantee which depends linearly on the neural network's prediction error, assuming relatively complete recourse.
    
    \item We perform a comprehensive experimental evaluation on two-stage robust knapsack, capital budgeting, and facility location benchmarks from the 2RO literature.~\method{} finds solutions of similar quality to or better than state-of-the-art. Large instances benefit the most from our method, with $100\times$ reduction and $5-10\times$ reductions in running time for knapsack and capital budgeting, respectively. 
    On the facility location problem, which involves objective and constraint uncertainty,
    ~\method{} computes higher quality solutions than $k$-adaptability in one to two orders of magnitude less time.  
    Moreover, for larger instances, $k$-adaptability is only tractable with $k=1$, which results in very conservative solutions, whereas \method{} achieves a rate of 89\% feasibility, which can be improved to 98\% by augmenting \method{} with 
    additional cuts. 
\end{enumerate}

The paper is organized as follows: \Cref{sec:background} provides the formal definitions of the general mixed-integer linear 2RO problem and the CCG algorithm. In~\Cref{sec:related_work}, we review the relevant literature on adjustable robust optimization, machine learning for robust optimization, and MILP representations of neural networks.~\Cref{sec:neur2ro} presents the~\method{} framework.~\Cref{sec:model} details the deep learning architecture of~\method{} and~\Cref{sec:theory} deals with its theoretical guarantees. The experimental results are presented in~\Cref{sec:experiments} before concluding in~\Cref{sec:conclusion}.

\section{Preliminaries}
\label{sec:background}
\subsection{Two-stage robust optimization}\label{sec:background_2RO}
Two-stage robust problems involve two types of decisions. The first set of decisions, $\vect x$, are referred to as \textit{here-and-now} decisions (or first-stage decisions) and are made before the uncertainty is realized. The second set of decisions, $\vect y$, are referred to as the \textit{wait-and-see} decisions (or second-stage decisions) and can be made on the fly after the uncertainty is realized. Note that this makes 2RO {less conservative} than standard RO due to the adjustability of the wait-and-see decisions. The uncertain parameters $\bxi$ are assumed to be contained in a convex and bounded uncertainty set $\Xi\subset \R^q$. The 2RO problem seeks a first-stage solution $\vect x$ which minimizes the worst-case objective value over all scenarios $\bxi\in \Xi$, where for each scenario the best possible second-stage decision $\vect y$ is implemented.

\paragraph{Example (Two-Stage Robust Facility Location).}
As a classical example of 2RO, consider the two-stage robust facility location problem. A planner must open a subset of $n$ facilities to serve $m$ clients with uncertain demands while minimizing the cost of opening the facilities and serving the clients. A facility with capacity $C_i$ can be opened at cost $c_i$. A customer has a demand $b_j(\bxi)$ that depends on an uncertain scenario $\bxi$.  Specifically, for a nominal demand $\bar{b}_j$ and some deviation $\Delta_j$, the demand is $b_j(\bxi) = \bar{b}_j + \Delta_j \cdot \xi_j$ where we assume there is a budget $B$ on the total percental amount of deviation, modeled by the budgeted uncertainty set~$\Xi = \{\bxi \in {[0, 1]}^{m} : \sum_{j=1}^{m} \xi_j \leq B\}$. Facility $i$ can serve customer $j$ at cost $d_{ij}$. The planner must decide which facilities to open before the uncertain demands are observed.  However, allocating facilities to the customers can be done \textit{after} the uncertain (demand) scenario is realized. This results in the following formulation:%
\begin{subequations}  \label{eq:FL}%
\begin{align}
    \min_{\vect x \in \gX} \max_{\bxi \in \Xi} \min_{\vect y \in \gY} \quad & \vect c^\intercal \vect x + \vect{d}^\intercal \vect y  & \\ 
    \text{s.t.} \quad  
        & \sum_{j=1}^m y_{ij} \cdot (\bar{b}_j + \Delta_j \cdot \xi_j)  \leq C_i \cdot x_i, & \forall i \in \{1,\ldots,n\},  \label{eq:FL_budget_constr}\\ 
        & \sum_{i=1}^n y_{ij} = 1, & \forall j \in \{1,\ldots,m\} \label{eq:FL_demand_constr}.
\end{align}
\end{subequations}
Here, $\gX = \{0, 1\}^n$ and entry $x_i$ of $\vect x$ is one if facility $i$ is opened in the first stage;  $\gY = \{0, 1\}^{n\times m}$ and entry $y_{ij}$ of $\vect y$ is one if facility $i$ is assigned to customer $j$ in the second-stage. 
Constraints~\eqref{eq:FL_budget_constr} and~\eqref{eq:FL_demand_constr} enforce that the facility capacities are not exceeded and ensure that each customer is served by exactly one facility, respectively.
The planner is concerned with finding the subset of facilities that minimizes the total cost over all scenarios while ensuring that the selected facilities can serve the customers for \textit{any} demand scenario, i.e., the first-stage decision $\vect x$ admits a feasible second-stage decision $\vect y$ for all $\bxi\in\Xi$.
In addition to demand uncertainty, Appendix~\ref{app:fl} presents a formulation that also considers uncertain facility disruptions on which we perform experiments. 

\paragraph{General 2RO Formulation.}
Mathematically, we study 2RO problems of the form%
\begin{subequations}
\begin{align}
    \min_{\vect x\in\gX} \max_{\bxi\in\Xi} \min_{\vect y\in\gY}\quad & 
        \vect c^\intercal \vect x + \vect d(\bxi)^\intercal \vect y \label{eq:2s_obj}\\
    \text{s.t.} \quad & T(\bxi)\vect x + W(\bxi) \vect y \le \vect h(\bxi), \label{eq:2s_cst}
\end{align}
\label{eq:2s_ro}%
\end{subequations}
where $\gX \subset \mathbb{Z}^n$ and $\gY \subset \mathbb{Z}^m$ are integer sets. The objective function and constraints of the problem, $\vect d(\bxi)\in \R^m, T(\bxi)\in\R^{r\times n}, W(\bxi)\in \R^{r\times m}$,  and $\vect h(\bxi)\in\R^r$, may depend on the scenario $\bxi$.
The first-stage objective coefficients, $\vect c\in\R^n$, may also be uncertain, but as \citet{subramanyam2020k} show, the uncertain parameters can be shifted to the second stage. Hence, without loss of generality, we assume the cost vector $\vect c$ to be deterministic.  
Since uncertainty in the objective function can be moved to the constraints by an epigraph reformulation, constraint uncertainty is the more general case. However, prior research on the special case of objective uncertainty has lead to more tractable algorithms. In our experiments, we study problems with objective and/or constraint uncertainty. Note that for constraint uncertainty, it can happen that for a given first-stage solution $\vect x$, there exists a scenario in $\Xi$ such that no feasible second-stage solution $\vect y$ exists. In this case, with slight abuse of notation, we assume that the inner minimization problem has an optimal value of $\infty$, enforcing that such solutions $\vect x\in\mathcal X$ cannot be optimal. We assume that there always exists a solution $\vect x\in\mathcal X$ which has finite objective value.

Our approach applies to the mixed-integer case as well but we will focus on the pure integer case here for ease of presentation and because all three benchmark problems are as such. 

\subsection{Column-and-Constraint Generation} \label{sec:ccg}

\let\oldnl\nl
\newcommand{\nonl}{\renewcommand{\nl}{\let\nl\oldnl}}

\begin{algorithm}[htbp!]
\footnotesize
\nonl \textbf{Input:} a problem instance's certain parameters;
 $\Xi$, the uncertainty set. \\
 \nonl \textbf{Output:} $\vect x^\star$, an optimal solution to \eqref{eq:2s_ro}.
\begin{algorithmic}[1]
\STATE set $ub=\infty$, $lb=-\infty$
\STATE $\Xi'=\{\bxi\}$ for any $\bxi\in \Xi$
 \WHILE{$ub-lb>0$}
    \STATE Calculate an optimal solution $(\vect x^\star,\mu^\star)$ to the main problem \eqref{eq:mp_levelset} and set $lb=\mu^*$. 
    \STATE Calculate an optimal solution $\bxi^\star$ (with optimal value $opt^*$) to the adversarial problem \eqref{eq:adversarial_appendix}.
    \STATE Set $\Xi'=\Xi'\cup \{\bxi^\star\}$ and $ub = \min\{ub, opt^\star\}$.
 \ENDWHILE
 \RETURN $\vect x^\star$
 \caption{Vanilla Column-and-Constraint Generation}
 \label{alg:CCG}
 \end{algorithmic}%
\end{algorithm}
CCG alternates between the \textit{main problem} and the \textit{adversarial problem}. We summarize its steps in~\Cref{alg:CCG}.
The main problem is given by%
\begin{subequations}
\begin{align*}
    \min_{\vect x\in\gX} \max_{\bxi\in\Xi'} \min_{y\in\gY} & 
        \quad \vect c^\intercal \vect x + \vect d (\bxi)^\intercal \vect y \nonumber\\
    \text{s.t.} & 
        \quad T(\bxi) \vect x + W(\bxi) \vect y  \le \vect h(\bxi), \nonumber&
\end{align*}
\end{subequations}
where $\Xi'\subset \Xi$ is a finite subset of scenarios. The main problem provides a lower bound for the optimal value of~\eqref{eq:2s_ro} as the maximum is taken only over a subset of scenarios. Using a level-set reformulation and introducing copies of the second-stage variables for each scenario in $\Xi'$, the problem can be written as%
\begin{subequations}
\begin{align}
    \min_{\vect x, \vect y^{(\bxi)}, \mu} \quad & \vect c^\intercal \vect x + \mu\\
    \text{s.t.} \quad & 
     \vect d(\bxi)^\intercal \vect y^{(\bxi)}\le \mu,   &\forall \bxi\in \Xi',\\
        & T(\bxi) \vect x + W(\bxi) \vect y^{(\bxi)} \le \vect h(\bxi),  &\forall \bxi\in\Xi',\label{eq:mp_levelset_cons}\\
        & \vect x\in\gX, \vect y^{(\bxi)}\in \gY,  &\forall \bxi\in\Xi'.
\end{align}
\label{eq:mp_levelset}%
\end{subequations}
In each iteration of CCG, a MILP solver computes an optimal solution $(\vect x^{\star},\mu^{\star})$ to \eqref{eq:mp_levelset}. The adversarial problem is solved next to find a new scenario that cuts off the current solution $\vect x^{\star}$:
\begin{subequations}
\begin{align}
     \bxi^\star\in\argmax_{\bxi\in\Xi} \min_{\vect y\in\gY} \quad & 
          \vect  d(\bxi) ^\intercal \vect y \\
    \text{s.t.} \quad & 
        W(\bxi) \vect y \le \vect h(\bxi) - T(\bxi) \vect x^\star . 
\end{align}
\label{eq:adversarial_appendix}%
\end{subequations}
When no such scenario can be found, the current solution $\vect x^{\star}$ is declared optimal. 

CCG often fails to converge quickly since both subproblems are hard to solve. In each iteration, the size of the former increases to accommodate new scenarios. The (bilevel) adversarial problem is extremely challenging for second-stage problems with integer variables. For mixed-integer second-stage problems,~\citet{zhao2012exact} proposed a nested CCG. Besides being computationally prohibitive, the nested CCG approach is not applicable to purely integer second-stage problems such as the ones we consider here.

\section{Related Work}
\label{sec:related_work}
\subsection{Algorithms for adjustable robust optimization}
Both single- and multi-stage robust mixed integer problems can be NP-hard even for problems that can be solved in polynomial time in its deterministic version \citep{buchheim2018robust}. Compared to single-stage robust problems, which are often computationally tractable as they can be solved using reformulations \citep{ben2009robust} or constraint generation \citep{mutapcic2009cutting}, two-stage problems are much harder to solve. When dealing with integer first-stage and continuous recourse, CCG is one of the key approaches \citep{zeng2013solving,tsang2023inexact}. However, many real-world applications involve integer second-stage decisions. While an extension of CCG has been proposed for mixed-integer recourse \citep{zhao2012exact}, it is often intractable and does not apply to pure integer second-stage problems. 

For problems with only objective uncertainty, 2RO can be solved by oracle-based branch-and-bound methods \citep{kammerling2020oracle}, branch-and-price \citep{arslan2022decomposition}, or iterative cut generation using Fenchel cuts \citep{detienne2024adjustable}. For special problem structures and binary uncertainty sets, a Lagrangian relaxation can be used to transform 2RO problems with constraint uncertainty into 2RO with objective uncertainty which can then be solved by the aforementioned methods \citep{subramanyam2022lagrangian,lefebvre2023adjustable}.

Besides these exact solution methods, several heuristics have been developed, including $k$-adaptability \citep{bertsimas2010finite,hanasusanto2015k,subramanyam2020k,ghahtarani2023double,arslan2022min,kurtz2024approximation}, decision rules \citep{bertsimas2018binary,bertsimas2015design}, and iteratively splitting the uncertainty set \citep{postek2016multistage}. 
In this paper, the $k$-adaptability branch-and-bound algorithm of \citet{subramanyam2020k} is used as a baseline as it is one of the only methods that can find high-quality solutions for reasonably large problems. 

\subsection{ML for robust optimization}
\citet{goerigk2022data} train a decision tree classifier to predict good start scenarios for CCG. Their method is orthogonal to ours in that it does not interfere with the iterations of CCG and so could be used in conjunction with~\method{}. However, it has only been applied to small shortest path and traveling salesperson problems on graphs with 6--9 nodes; this is because data collection requires solving the training instances to optimality, a strategy that cannot be scaled. 
\citet{julien2024machine} use ML to improve node selection in the $k$-adaptability branch-and-bound algorithm of \citet{subramanyam2020k}. For capital budgeting, slightly better solutions are found earlier on in the tree search compared to the vanilla version of $k$-adaptability. This method also requires solving training instances using a randomized branch-and-bound algorithm, a limiting factor for scaling up. In comparison,~\method{}'s data collection reduces to solving (typically ``easy'') MILPs for the second-stage problem which can be parallelized across samples.

\citet{bertsimas2024machine} use decision tree classifiers to predict first-stage decisions as a function of the parameters of an instance. They only consider problems with continuous second-stage variables whereas~\method{} can handle integers in both stages. Similar to our setting, training instances are required for training. Contrary to our method,~\citet{bertsimas2024machine} need optimal or near-optimal solutions to the training instances (e.g., using CCG) to fit the decision tree(s).  Additionally, as the decision tree can only output a limited number of first-stage decisions, these have to be sufficiently representative of the set of all optimal solutions; problems with a high sensitivity to instance parameters would be difficult to handle. Another key limitation is that some problem parameters are assumed to be fixed across instances. In the facility location problem, for example, this means that all training and test instances must have the same number of facilities $n$ and customers $m$ and that predictions are not invariant to the arbitrary ordering of the facilities/customers. Our neural network architecture does not suffer from this limitation. 

Besides improving algorithmic performance, ML has been used to construct uncertainty sets based on historical data~\citep{goerigk2023data,tulabandhula2014robust,ning2018data,shang2017data,shang2019data,shen2020large,ning2017data,campbell2015bayesian,wang2023learning}. 
While interesting and related, here we assume that the uncertainty set is known.

\subsection{MILP representations of neural networks}
One key aspect of \method{} is representing a trained neural network (with rectified linear units as activations) as MILPs, an idea first explored in \citet{cheng2017maximum, tjeng2017evaluating, fischetti2018deep}. These representations are being improved~\citep{grimstad2019relu, anderson2020strong, wang2023optimizing} and packaged into open software \citep{bergman2022janos,ceccon2022omlt,turner2023pyscipopt}, including by Gurobi~\citep{gurobiml}. Neural network surrogates have been proposed for non-linear or intractable constraints \citep{say2017nonlinear,grimstad2019relu,murzakhanov2020neural,KATZ202011350,KODY2022108282,kotary2021learning,chen2022learning}, stochastic programming \citep{dumouchelle2022neur2sp,alcantara2024quantile,wuhgcn2sp,kronqvist2023alternating,bae2023deep}, and bilevel optimization \citep{dumouchelle2024neur2bilo,zhou2024learning}.  
The min-max-min optimization that~\method{} deals with has not been explored in these previous works, necessitating deep integration with CCG, as well as several neural network design choices specific to 2RO.

\section{\method{}: A Learning-augmented CCG Algorithm}
\label{sec:neur2ro}

\subsection{Problem setting}
Our approach aims to train {neural networks} that estimate the optimal value and the degree of infeasibility of the second-stage problem of \eqref{eq:2s_ro} for a fixed first-stage decision and scenario. 
In this section, we explicitly define two networks as it is the most general view.  However, in Section~\ref{sec:model}, we present a single architecture that makes two predictions using almost entirely shared parameters.
The {neural networks} are then integrated within CCG {to obtain first-stage decisions that have good objective values as predicted by the value network and are likely second-stage-feasible as predicted by the feasibility network}. We rely on a training dataset of historical instances, as is typically assumed in ML-for-optimization works. More explicitly, consider a~\textit{parametric} 2RO setting:%
\begin{subequations}
\begin{align}
    \min_{\vect x\in\gX} \max_{\bxi\in\Xi} \min_{\vect y\in\gY}\quad & 
        \vect c (\bpi)^\intercal \vect x + \vect d(\bxi,\bpi)^\intercal \vect y \label{eq:2s_obj_parametric}\\
    \text{s.t.} \quad & T(\bxi,\bpi)\vect x + W(\bxi,\bpi) \vect y \le \vect h(\bxi,\bpi). \label{eq:2s_cst_parametric}
\end{align}
\label{eq:2s_ro_parametric}%
\end{subequations}
Compared to Problem~\eqref{eq:2s_ro}, Problem~\eqref{eq:2s_ro_parametric} makes explicit the dependence of the objective function and constraints on~\textit{certain, observed} parameters $\bpi\in\Pi$.
Generally, $\gX$, $\gY$, and $\Xi$ depend on $\bpi$ as we consider varying uncertainty budgets and instance sizes.  However, for brevity, we omit this dependence in the notation.
In the facility location problem, $\bpi$ refers to the nominal demands $\bar{b}_j$, deviation factors $\Delta_j$ of each customer, capacities $C_i$ of each facility, facility opening costs $\vect c$, serving costs $\vect d$, number of facilities, number of customers, and the uncertainty budget.

Instances of the problem are assumed to be drawn i.i.d. from a fixed but unknown distribution over the set of certain parameters $\Pi$.  Moreover, the uncertainty set $\Xi$ is known. At training time,~\method{} has access to a sample of instances from that distribution. At test time, unseen instances from the same distribution are used for evaluation. 

\subsection{Neural networks for objective and feasibility prediction}

\method{} attempts to replace the main~\eqref{eq:mp_levelset} and adversarial~\eqref{eq:adversarial_appendix} problems with surrogates that are easier to solve. It does so using two MILP-representable neural networks with learnable parameters $\Theta$ that can assess, for a first-stage decision $\vect{x}$ under uncertain scenario $\bxi$ and certain instance parameters $\bpi$, the following:

\begin{enumerate}[label=(\alph*)]
    \item {$\NN_\Theta^F(\bpi, \vect x, \bxi)$} predicts whether there is a feasible solution $\vect{y}$ to the second-stage problem parameterized by $\vect{x}$ and $\bxi$. It does so by predicting the minimum amount of slack $\beta$ needed to render the constraints feasible, i.e., 
    \[\NN_\Theta^F(\bpi, \vect x, \bxi) \approx\min_{\vect y\in\gY, \beta\in\R_{+}}\{\beta : \boldsymbol{1}^\intercal\beta \ge \vect h(\bxi, \bpi) - W(\bxi, \bpi) \vect y - T(\bxi, \bpi) \vect x\}.\] 
    Note that $\beta=0$ if $\vect{x}$ is feasible;
    \item {$\NN_\Theta^O(\bpi, \vect x, \bxi)$} predicts the second-stage objective value of (a feasible) $\vect{x}$, i.e., \[\NN_\Theta^O(\bpi, \vect x, \bxi) \approx\min_{\vect y\in\gY}\{\vect d(\bxi,\bpi)^\intercal \vect y : W(\bxi,\bpi) \vect y \le \vect h(\bxi,\bpi) - T(\bxi,\bpi) \vect x \}.\] 
\end{enumerate} 
Why are two predictions required? If we consider that the ``value'' of an infeasible second-stage (minimization) problem is $+\infty$, then one could train a single value prediction network with output range in $\mathbb{R}\cup\{+\infty\}$, where $\infty$ could be represented by a big-$M$ value. However, this would lead to fitting a discontinuous function by a continuous neural network function, leading to large errors on the feasibility boundary. Since the cost of mispredicting $M$ is extremely large, this would lead to a model that is likely to predict most inputs to be infeasible. To resolve this issue, we use one model for each of the two tasks. For problems with relatively complete recourse, the feasibility prediction network is not used.

Assuming already trained neural networks with parameters $\Theta^\star$ that produce accurate predictions for the latter two tasks, the \textit{surrogate} main problem can be formulated as%
\begin{subequations}
\begin{align}
    \min_{\vect x\in\gX, \mu} 
       \quad  & c (\bpi)^\intercal \vect x + \mu \\
    \text{s.t.} \quad &   \NN^O_{\Theta^{\star}}(\bpi,\vect x, \bxi) \leq\mu,   &\forall \bxi\in \Xi', \label{eq:max_nn_obj}\\
    &  \NN^F_{\Theta^{\star}}(\bpi,\vect x, \bxi) \leq\varepsilon, &\forall \bxi\in \Xi', \label{eq:max_nn_feas}
\end{align}
\label{eq:max_nn}%
\end{subequations}
where $\varepsilon\in\mathbb{R}_{>0}$ is a small user-defined feasibility tolerance. Notice that the per-scenario second-stage variables $\vect{y}^{(\bxi)}$ are not part of this surrogate, a key to reducing the complexity of the main problem as we will discuss later. 
This surrogate problem seeks a first-stage solution $\vect x$ that satisfies two conditions: it minimizes the maximum (i.e., worst-case) output of the objective neural network across all scenarios in $\Xi'$ (constraint~\eqref{eq:max_nn_obj}); it is predicted to be $\varepsilon$-feasible by the feasibility neural network across all scenarios in $\Xi'$ (constraint~\eqref{eq:max_nn_feas}). If both networks are perfectly accurate, then an optimal solution to~\eqref{eq:max_nn} is also optimal in the original main problem~\eqref{eq:mp_levelset} if we choose $\varepsilon=0$. 

The corresponding surrogate model of the adversarial problem is no longer an intractable bilevel problem as in~\eqref{eq:adversarial_appendix}. It simply requires optimizing over the inputs of the two neural networks: 
\begin{equation}\bxi^{(a)} \in \argmax_{\bxi \in \Xi} \; \NN^O_{\Theta^{\star}}(\bpi, \vect x^\star, \bxi), \quad \bxi^{(b)} \in \argmax_{\bxi \in \Xi} \; \NN^F_{\Theta^{\star}}(\bpi, \vect x^\star, \bxi).\label{eq:ap_o_}\end{equation}

\subsection{Data collection and model training}
\Cref{alg:training} summarizes the data collection and model training processes for the two neural networks.

\begin{algorithm}[htbp!]
\footnotesize
\nonl \textbf{Input:} $\Pi_{\text{train}}=\{\bpi^{(k)}\}\subset\Pi$, a set of sampled or historical parameter vectors; $\Xi$, the uncertainty set.\\
\nonl \textbf{Output:} $\Theta^{\star}_O$, $\Theta^{\star}_F$, neural network parameters for objective/feasibility prediction.
\begin{algorithmic}[1]
 \STATE Initialize $\mathcal{D}_{\text{train}}^O$ and $\mathcal{D}_{\text{train}}^F$ to empty sets
 \FORALL{training instances $\bpi^{(k)}\in\Pi_{\text{train}}$}
    \STATE $\mathcal{X}_\text{train}=\left\{ \vect{x}^{(i)} \overset{U}{\sim}\mathcal{X}\right\}_{i=1}^{T}$, a set of $T$ first-stage decision vectors drawn from $\mathcal{X}$ 
    \FORALL{$\vect{x}^{(i)} \in \mathcal{X}_\text{train} $}
        \STATE $\Xi_\text{train}=\left\{ \bxi^{(j)} \overset{U}{\sim}\Xi\right\}_{j=1}^{L}$, a set of $L$ scenarios drawn from $\Xi$
        \FORALL{pairs $(\vect{x}^{(i)}, \bxi^{(j)})\in \mathcal{X}_\text{train} \times \Xi_\text{train}$}
            \STATE Solve the second-stage problem:
            $$\min_{\vect y\in\gY}\{ \vect d(\bxi^{(j)},\bpi^{(k)})^\intercal \vect y : W(\bxi^{(j)},\bpi^{(k)}) \vect y \le \vect h(\bxi^{(j)},\bpi^{(k)}) - T(\bxi^{(j)},\bpi^{(k)}) \vect x^{(i)} \}$$
            \STATE If feasible, let $\lambda^{(i,j)}$ be its optimal value and set $\beta^{(i,j)}=0$.
            \STATE If infeasible, solve the following problem and let $\beta^{(i,j)}$ be its optimal value:
            $$\min_{\vect y\in\gY, \beta\in\mathbb{R}_+}\{\beta : \boldsymbol{1}^\intercal\beta \ge \vect h(\bxi^{(j)}, \bpi^{(k)}) - W(\bxi^{(j)}, \bpi^{(k)}) \vect y - T(\bxi^{(j)}, \bpi^{(k)}) \vect x^{(i)}\}.$$
        \ENDFOR
        \STATE For all pairs $i,j$ with $\beta^{(i,j)} = 0$, add $((\bpi^{(k)},\vect{x}^{(i)},\bxi^{(j)}),\lambda^{(i,j)})$ to $\mathcal{D}_{\text{train}}^O$ 
        \STATE For all pairs $i,j$ add $((\bpi^{(k)},\vect{x}^{(i)},\bxi^{(j)}),\beta^{(i,j)})$ to $\mathcal{D}_{\text{train}}^F$ 
    \ENDFOR
 \ENDFOR
 \STATE Find the best-fit network parameters for objective prediction by minimizing the mean squared error: 
 $$\Theta^{\star}_O = \argmin_{\Theta}\sum_{\left(\left(\bpi^{(k)},\vect{x}^{(i)},\bxi^{(j)}\right),\lambda^{(i,j)}\right)\in\mathcal{D}_{\text{train}}^O}(\lambda^{(i,j)} - \NN^O_\Theta(\bpi^{(k)},\vect{x}^{(i)},\bxi^{(j)}))^2,$$
 \STATE and similarly for feasibility prediction: 
 $$\Theta^{\star}_F = \argmin_{\Theta}\sum_{\left(\left(\bpi^{(k)},\vect{x}^{(i)},\bxi^{(j)}\right),\beta^{(i,j)}\right)\in\mathcal{D}_{\text{train}}^F}(\beta^{(i,j)} - \NN^F_\Theta(\bpi^{(k)},\vect{x}^{(i)},\bxi^{(j)}))^2$$ 
 \RETURN $\Theta^{\star}_O$ and $\Theta^{\star}_F$
 \caption{Data collection and model training}
 \label{alg:training}
 \end{algorithmic}
\end{algorithm}

\subsubsection{Sampling training data}\label{sec:sampling}
Given an instance $\bpi$, we would like to sample first-stage decisions $\vect x\in\mathcal{X}$ and scenarios $\bxi\in\Xi$ that provide good coverage over the range of possible second-stage outcomes. For 2RO problems with only objective uncertainty, this translates into good coverage of the range of second-stage objective values. When there is uncertainty in the constraints, we additionally would like to sample first-stage decisions and scenarios for which the second-stage problem is infeasible so that we can learn to avoid such first-stage decisions.  

We will assume that the first-stage variables are binary as this is the case for many applications. This is without loss of generality since our method also applies to the mixed-integer case. 
To sample binary first-stage decisions, we first draw a value $p\sim U(0,1)$, then independently set each variable $x_i$ to one with probability $(1-p)$ and to zero otherwise.

To sample scenarios from box-constrained uncertainty sets, e.g., ~$\Xi = \{\bxi \in {[0, 1]}^{q} \}$, we draw a value for each dimension $\xi_i$ of the scenario vector independently and uniformly at random from the bounded domain of the scenario.
For budgeted uncertainty sets with, e.g., $\Xi = \{\bxi \in {[0, 1]}^{q} : \sum_{j=1}^{q} \xi_j \leq B\}$, such as in facility location, we sample $\xi_i \sim U(0,1)$ then normalize the vector $\bxi$ to sum up to the uncertainty budget $B$, i.e., $\bxi = B \cdot \bxi / \sum_{i=1}^q \xi_i$. This guarantees that scenarios in the training set are sufficiently difficult to satisfy in the second-stage problem. Additional problem-specific sampling measures are deferred to~\Cref{app:sampling}.

For each instance $\bpi^{(k)}$ in the training set $\Pi_{\text{train}}$, $T$  first-stage decisions and $L$ scenarios per first-stage decision are sampled as we just described (lines 1--5 of~\Cref{alg:training}). 
Then, for each pair, the second-stage problem is solved and its objective value (if feasible) or amount of constraint violation (if infeasible) are recorded as targets for supervised learning (lines 6--10). 

\subsubsection{Model training} 
In the general case of constraint uncertainty, both neural networks are trained by minimizing a regression loss function (here, mean squared error, lines 9--10). The resulting neural networks are simply real-valued functions that we denote by $\NN^O_{\Theta^{\star}}(\cdot)$ and  $\NN^F_{\Theta^{\star}}(\cdot)$. In~\Cref{sec:model}, we will discuss the neural network architecture in detail. Since the two prediction tasks share the same inputs and only differ in their target outputs, they will share most of the parameters, resulting in basically a single two-output model.

\subsection{Integration of trained models into CCG} 

\begin{algorithm}[htbp!]
\footnotesize
\nonl\textbf{Input:} $\bpi$, a problem instance's certain parameters;
$\Xi$, the uncertainty set; 
$\Theta^{\star}_O$, $\Theta^{\star}_F$, neural network parameters for objective/feasibility prediction; $\varepsilon>0$, an accuracy parameter.\\
\nonl\textbf{Output:} $\vect x^\star \in \gX$.%
\begin{algorithmic}[1]
\STATE $\converged = \false$
\STATE $\Xi_F'=\Xi_O'=\{\bxi\}$ for any $\bxi\in \Xi$ 
 \WHILE{not \converged}
        \STATE Calculate an optimal solution $\vect{x}^{\star}$ 
        of the main problem \eqref{eq:max_nn2}
        \STATE Calculate optimal solutions $\bxi_O$ and $\bxi_F$ of the adversarial problems
        \begin{minipage}{.95\linewidth}
               \begin{equation}\bxi_O \in \argmax_{\bxi \in \Xi} \; \NN^O_{\Theta^{\star}}(\bpi, \vect x^\star, \bxi),   \label{eq:ap_o} \end{equation}        
               \begin{equation}\bxi_F \in \argmax_{\bxi \in \Xi} \; \NN^F_{\Theta^{\star}}(\bpi, \vect x^\star, \bxi). \label{eq:ap_f}\end{equation}
        \end{minipage}

        \STATE $\converged = \true$
        \IF{$\NN^O_{\Theta^{\star}}(\bpi, \vect x^\star, \bxi_O) \ge \max_{\bxi\in \Xi_O'} \NN^O_{\Theta^{\star}}(\bpi,\vect x^\star, \bxi) + \varepsilon$}
            \STATE Update $\Xi_O'=\Xi_O'\cup \{\bxi_O\}$
            \STATE $\converged = \false$
        \ENDIF

        \IF{$\NN^F_{\Theta^{\star}}(\bpi, \vect x^\star, \bxi_F) \ge \varepsilon$}
            \STATE Update $\Xi_F'=\Xi_F'\cup \{\bxi_F\}$
            \STATE $\converged = \false$
        \ENDIF
    \ENDWHILE
 \RETURN $\vect x^\star$
 \caption{\method{}'s Column-and-Constraint Generation}
 \label{alg:CCG_neur2ro}
 \end{algorithmic}
\end{algorithm}
\Cref{alg:CCG_neur2ro} shows how the trained neural networks are incorporated into CCG; the original CCG is shown in~\Cref{alg:CCG} for comparison.~\Cref{alg:CCG_neur2ro} maintains two sets of scenarios, $\Xi_O'$ and $\Xi_F'$, that are grown with new violating scenarios at each iteration; they are initialized with an arbitrary scenario in Step 2. The algorithm then follows the structure of the vanilla CCG. 

\subsubsection{Main problem} 
We have already considered an initial version of the surrogate main problem in~\eqref{eq:max_nn}.

One drawback of this surrogate is that it is very sensitive to the accuracy of the neural networks, since,  even if the considered subset of scenarios $\Xi'$ is the optimal choice, the optimal value and feasibility constraints highly depend on the predicted (and probably slightly inaccurate) values of the neural networks. To alleviate this issue, we developed a new formulation where the neural networks only decide on the most critical scenarios for the objective function and the constraints rather than the exact objective/infeasibility estimates by using an $\argmax$ formulation. The predicted worst-case scenarios are then used to calculate the exact objective value and to enforce feasibility. In this approach, the neural networks only have to be accurate in ``ranking'' the given scenarios; slight inaccuracies in the predicted value do not change the solution of the problem as long as the ranking remains correct. This idea is incorporated in the following formulation.
\begin{subequations}
\begin{align}
    \min_{\substack{\vect x\in\gX, \mu,\\\bxi^{(a)}\in\Xi_O', \bxi^{(b)}\in\Xi_F'}} 
       \quad  & \mu \\
    \text{s.t.} \quad & \mu=\vect c (\bpi)^\intercal \vect x + \vect d(\bxi^{(a)},\bpi)^\intercal \vect y, \label{eq:max_nn_obj_}\\
    & W(\bxi^{(a)},\bpi) \vect y + T(\bxi^{(a)},\bpi) \vect x \le \vect h(\bxi^{(a)},\bpi), \label{eq:mp_objfeas_neur2ro_}\\ 
    & W(\bxi^{(b)},\bpi) \vect y + T(\bxi^{(b)},\bpi) \vect x \le \vect h(\bxi^{(b)},\bpi), \label{eq:mp_consfeas_neur2ro_}\\ 
    &  \bxi^{(a)} \in {\arg\max}_{\bxi \in \Xi_O'} \big\{\NN^O_{\Theta^{\star}}(\bpi,\vect x, \bxi)  \big\}\label{eq:mp_argmax1_neur2ro_}, \\
    &  \bxi^{(b)} \in {\arg\max}_{\bxi \in \Xi_F'} \big\{\NN^F_{\Theta^{\star}}(\bpi,\vect x, \bxi)  \big\}.\label{eq:mp_argmax2_neur2ro_}
\end{align}
\label{eq:max_nn2}%
\end{subequations}
Constraints \eqref{eq:mp_argmax1_neur2ro_} and \eqref{eq:mp_argmax2_neur2ro_} select the predicted worst-case scenarios $\bxi^{(a)}$ and $\bxi^{(b)}$ for the objective value and the feasibility, respectively. Both scenarios are plugged into the constraint system in \eqref{eq:mp_objfeas_neur2ro_} and \eqref{eq:mp_consfeas_neur2ro_} to ensure feasibility in both scenarios. In \eqref{eq:max_nn_obj_} the objective value is evaluated on the worst-case scenario $\bxi^{(a)}$. In contrast to the classical main problem \eqref{eq:mp_levelset}, we only need a single second-stage vector $\vect y$, which is independent of the number of scenarios added to the problem. Note that the $\argmax$ constraints can be modeled by mixed-integer linear constraints involving big-$M$ constraints; see Section \ref{sec:bigm}. We empirically demonstrate the efficacy of \eqref{eq:max_nn2} over \eqref{eq:max_nn} in~\Cref{app:ab_argmax}.

\subsubsection{Adversarial problem} 
With a candidate first-stage solution $\vect{x}^{\star}$ in hand, we can move on to the adversarial problems in Step 5. In contrast to vanilla CCG, we find up to~\textit{two} violating scenarios for $\vect{x}^{\star}$ by solving one adversarial problem for each of the two neural networks. 
This is necessary because $\vect{x}^{\star}$ can be deficient in one of two ways: 
\begin{enumerate}[label=(\alph*)]
\item It is infeasible, i.e., there exists a scenario in $\Xi$ for which there does not exist a $\vect y\in\mathcal{Y}$ that satisfies the constraints with $\vect{x}^{\star}$. This is what the feasibility neural network $\NN^F_{\Theta^{\star}}$ is meant to detect;
\item It is feasible but its objective value is worse than currently evaluated, i.e., there exists a scenario in $\Xi$ for which the value of the second-stage problem is worse (larger) than that of any scenario in $\Xi_O'$. This is what the value network predicts.
\end{enumerate}
The two scenarios can be interpreted as feasibility and optimality cuts.
Problems~\eqref{eq:ap_o} and~\eqref{eq:ap_f} are mixed-integer programs (linear for polyhedral uncertainty sets) in which inputs $\bpi$ and $\vect x^\star$ to the neural networks are fixed and only the scenario vector $\bxi$ is variable. A violating scenario $\bxi_O$ is added to the set $\Xi_O'$ if the predicted second-stage objective value of the current iterate $\vect x^\star$ is~\textit{worse} than that of any scenario in $\Xi_O'$ (lines 7, 8). A scenario $\bxi_F$ is added to the set $\Xi_F'$ if the current iterate $\vect x^\star$ is predicted to be infeasible for the second-stage problem (lines 11, 12). If no scenario is added in an iteration, then the variable \converged{} triggers the termination of the algorithm. The last first-stage solution is returned.

This concludes the high-level description of our method. Compared to vanilla CCG~\PCref{alg:CCG}, the surrogate main problem is smaller in size as it does not require per-scenario copies of the instance's constraints and second-stage variables, but only of the neural networks; we will show in~\Cref{sec:model} that these can be encoded quite compactly. The surrogate adversarial problems are single-level mixed-integer programs in contrast to the typical bilevel formulation~\eqref{eq:adversarial_appendix}.

\section{A Neural Network Architecture for~\method}
\label{sec:model}
Next, we will describe the neural network architecture for objective predictions. Feasibility will ultimately be predicted using the same network, modulo a small tweak to allow for two outputs. 
\subsection{Required model properties}
For the predictive models, we are aiming for the following three properties:
\begin{enumerate}[label=(\alph*)]
    \item\textit{Predictive yet 
    efficient to optimize over:} 
    To achieve accurate predictions of the second-stage objective value, the neural network should have a certain expressivity, being able to approximate highly non-linear functions, which usually comes with a large number of learnable parameters. However, a more complex network structure leads to less tractable MILP representations in the main and adversarial problems of CCG (\Cref{alg:CCG_neur2ro}). Recall that the motivation for constructing surrogates of the main and adversarial problems is to get around the computational cost incurred in every iteration of CCG.
    \item\textit{Modularity:} 
    As per Step 4 of \Cref{alg:CCG_neur2ro}, multiple copies of the (trained) neural network---one for each scenario---will be encoded as a MILP in the surrogate main problem.  For this reason, the network architecture should be modular to allow only select components to be encoded rather than the entire network.   This enables scalability as fewer variables are introduced in each iteration of \Cref{alg:CCG_neur2ro}.
    \item\textit{Permutation and size invariance:} The network's prediction should be invariant to permutations of the decision variables, e.g., a re-indexing of the facilities and customers for facility location. It should also handle instances of different sizes, since, e.g., the number of customers in the facility location problem may change over time.
\end{enumerate}
Recall that the inputs of the predictive model are $\bpi$, $\vect x$, and $\bxi$. In the following we explore some basic machine learning models and see why they cannot satisfy the requirements (a)-(c).

First, a linear regression model is not able to model the highly non-linear dependencies between input $\bpi$, $\vect x$, $\bxi$, and the optimal value or constraint violation value. Note that even in the simplest case, when the second-stage problem is an LP, the optimal value function is piecewise-linear and convex in the right-hand-side constraint parameters and piecewise-linear and concave in the objective coefficients. Regarding the constraint parameters, the optimal value function may be even discontinuous. When considering integer decision variables, the functions to be fitted can be even less structured. Hence, linear models clearly violate (a).

Tree-based models, either single regression trees, random forests, or boosted tree ensembles, are MILP-representable~\citep{bertsimas2017optimal} and somewhat more interpretable than neural networks. While such tree-based models could in principle satisfy (a), they fail on (b) and (c). If we view such a model as a function $f(\bpi, \vect x, \bxi)$ where each of the three inputs is flattened into a vector of fixed dimensionality, then the size of $\vect x$ and the number of parameters in $\bpi$ and $\bxi$ must also be fixed, violating (c). Additionally, any re-indexing of the first-stage decision variables ``confuses'' the decision tree(s). A recent set-based decision tree was proposed by~\citet{hirsch2021trees} and could be explored in future work.

Alternatively, consider the simplest type of neural network, a fully-connected multilayer perceptron (MLP) that takes as input a concatenation of $\bpi, \vect x$, and $\bxi$ and predicts as $f(\text{concat}[\bpi, \vect x, \bxi])$. Clearly, this model violates (c); any permutation of the first-stage variable can change the output of the model. It also violates the modularity requirement in (b). Consider~\eqref{eq:mp_argmax1_neur2ro_} in~\Cref{alg:CCG_neur2ro} where the model predicts the value of first-stage variables $\vect x$ for a~\textit{fixed set of scenarios} $\Xi_O'$. The whole MLP has to be encoded as a MILP for each scenario, although the scenarios are fixed, not variable. Because this MLP will likely need a number of neurons proportional to the size of the input vector in each hidden layer to accurately predict the second-stage objective values (thus violating requirement (a)), encoding it many times in a MILP is not viable.

\subsection{A modular architecture}
The basic unit we will use here is an MLP with Rectified Linear Unit (ReLU) activation functions. We refer to Section 4.2.1 of~\citet{huchette2023deep} for MILP formulations of a trained ReLU MLP, including the most basic big-$M$ formulation we use in our implementation. We let $\Phi:\mathbb{R}^{q_\text{in}} \mapsto \mathbb{R}^{q_\text{out}}$
denote a ReLU MLP that maps a $q_\text{in}$-dimensional column vector to $\mathbb{R}^{q_\text{out}}$. An MLP has learnable parameters that will be optimized during training~\PCref{alg:CCG_neur2ro}. 

\begin{figure}[htbp!]
    \centering
    \includegraphics[width=\textwidth]{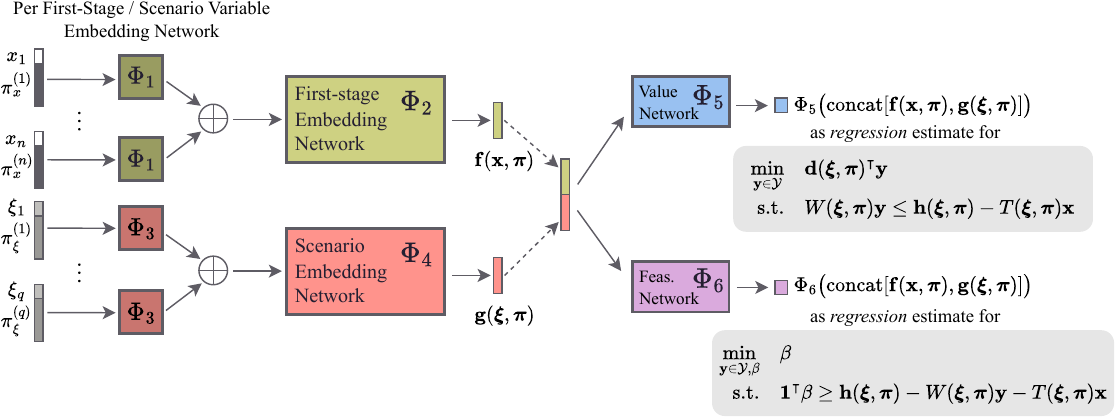}
    \caption{The neural network architecture for learning-augmented CCG. To distinguish between instance parameters related to the first-stage decisions and the scenario, we use $\pi_{\vect x}^{(i)}$ and $\pi_{\bxi}^{(j)}$ as parameters specific to the $i$-th first-stage decision and the $j$-th scenario index, respectively.  $\bigoplus$ denotes the aggregation (summation) of vectors.}
    \label{fig:full_architecture}
\end{figure}

The neural network architecture is illustrated in~\Cref{fig:full_architecture}.
The first part of the network structure is an embedding network that maps the instance parameters and a first-stage decision into a fixed-dimensional~\textit{embedding} (vector) such that 
\begin{enumerate*}[label=(\alph*)]
    \item a re-indexing of the variables does not change the embedding (permutation invariance), and
    \item first-stage problems with a varying number of variables can be tackled (size ``invariance'', a term we use loosely here to describe this property).
\end{enumerate*} 
Let $\vect x \in \R^n$ be a first-stage decision vector with entries $x_i$ and 
$\bpi_i$ a vector of certain parameters that are specific to the $i$-th variable (e.g., objective function or constraint coefficients of variable $x_i$) or generic to the instance (e.g., a right-hand side constraint coefficient). We concatenate these inputs into a single~\textit{feature vector} for each variable, $X^{(i)}=\text{concat}[x_i, \bpi_i]$. Each such feature vector is passed, independently of the other variables, through the same first embedding MLP: $$X^{\prime(i)} = \Phi_{1}\left(X^{(i)}\right)\in\mathbb{R}^{\text{dim}_1}.$$
We now have a set of vectors $\{X^{\prime(i)}\}_{i\in[n]}$ each of dimension $\text{dim}_1$. Each vector in this set is a fixed-size representation of a given variable $x_i$. 

We now need to capture interactions between the variables. As a first step, we ``pool'' these vectors together by aggregating them using a set function. For example, we can sum the vectors elementwise to obtain a single vector of dimension $\text{dim}_1$. Mathematically, this pooling operation is:
$$X^{\text{pool}}= \sum_{i=1}^n {X^{\prime(i)}}\in\mathbb{R}^{\text{dim}_1}.$$
Last, $X^{\text{pool}}$ is passed to another MLP $\Phi_{2}: \mathbb{R}^{\text{dim}_1} \mapsto \mathbb{R}^{\text{dim}_2}$ to increase predictive capacity:
$$X^{\text{final}} = \Phi_{2}\left(X^{\text{pool}}\right)\in\mathbb{R}^{\text{dim}_2}.$$

The composition of functions (MLPs and pooling) that maps the set of input vectors $\{X^{(i)}\}_{i\in[n]}$ to the single embedding vector  $X^{\text{final}}$ is permutation-invariant and can accommodate any number of first-stage variables $n$, as desired. The combination of a single shared MLP followed by pooling has origins in the DeepSets work of~\citet{zaheer2017deep}. Let us use $\vect{f} \left( \vect{x}, \bpi; \Phi_{1}, \Phi_{2}\right)$ to refer to the full composition which is parameterized by the two neural networks $\Phi_{1}(\cdot)$ and $\Phi_{2}(\cdot)$, i.e., $X^{\text{final}} = \vect{f} \left( \vect{x}, \bpi; \Phi_{1}, \Phi_{2}\right)$. 

The very same design can be used to embed a scenario $\bxi$. Each entry $\bxi_j$ of a scenario is represented by a feature vector $Z^{(j)}$ that depends on $\bxi$ and the certain parameters $\bpi$. A composition of MLPs, $\vect{g} \left( \bxi, \bpi; \Phi_{3}, \Phi_{4}\right)$, maps that set of feature vectors to an embedding of dimension $\text{dim}_4$, where $\Phi_{3}(\cdot)$ and $\Phi_{4}(\cdot)$ are the MLPs that map each scenario feature vector to $\mathbb{R}^{\text{dim}_3}$ and the result of their pooling to $\mathbb{R}^{\text{dim}_4}$, respectively. 

Finally, the two embeddings of the first-stage decision and the scenario are concatenated and passed on to an MLP $\Phi_{5}: \mathbb{R}^{\text{dim}_2+\text{dim}_4}\mapsto\mathbb{R}$ that predicts the objective value:
$$\NN^O_{\Theta}(\bpi,\vect x, \bxi) = \Phi_{5}\left(\text{concat}\left[\vect{f} \left( \vect{x}, \bpi; \Phi_{1}, \Phi_{2}\right), \vect{g} \left( \bxi, \bpi; \Phi_{3}, \Phi_{4}\right)\right]\right).$$

We simplify notation and reduce the computational cost by allowing the objective and feasibility predictors to share a common backbone. There is no need for completely separate predictors because both of them have the same inputs, namely an instance's certain parameters, a first-stage decision, and a scenario. We exploit this connection by letting both predictions be based on the same decision and scenario embedding vectors, giving a feasibility prediction function of the form: 
$$\NN^F_{\Theta}(\bpi,\vect x, \bxi) = \Phi_{6}\left(\text{concat}\left[\vect{f} \left( \vect{x}, \bpi; \Phi_{1}, \Phi_{2}\right), \vect{g} \left( \bxi, \bpi; \Phi_{3}, \Phi_{4}\right)\right]\right),$$
where $\Phi_{6}: \mathbb{R}^{\text{dim}_2+\text{dim}_4}\mapsto\mathbb{R}$ is another MLP. The two prediction functions differ only in the learnable parameters of their last MLPs, $\Phi_{5}$ and $\Phi_{6}$. Separate parameters allow for calibration of the predictions depending on the target. 

The subscript $\Theta$ to $\NN^F_{\Theta}$ and $\NN^O_{\Theta}$ refers to a set of values assigned to the parameters of MLPs $\Phi_1,\dots,\Phi_6$; $\Theta^\star$ are the parameter values obtained after training the whole network architecture using the usual stochastic gradient descent approach on the dataset collected in~\Cref{alg:training}.

\subsection{Optimizing over the neural network}
\label{sec:mp}

As mentioned earlier, a ReLU MLP with fixed parameters can be encoded as a MILP~\citep{huchette2023deep} which we use in the main and adversarial problems in~\Cref{alg:CCG_neur2ro}. Due to their modular structure, only some components of those networks need to be encoded in the main or adversarial problem. These savings in the number of additional decision variables and constraints make the overhead of optimizing over deep neural networks negligible relative to the original CCG algorithm.

\subsubsection{Main problem}
Specifically, in the relevant constraints~\eqref{eq:mp_argmax1_neur2ro_} and~\eqref{eq:mp_argmax2_neur2ro_} of the surrogate main problem, a prediction of the second-stage objective value is required for each $\bxi \in \Xi_O'$, and another prediction of feasibility for each $\bxi \in \Xi_F'$, respectively. To that end, for each scenario, a duplication of the MILP encoding of $\NN^{O}_{\Theta^\star}$ and $\NN^{F}_{\Theta^\star}$ has to be incorporated into the optimization formulation. Due to the modularity of our architecture and since the only variable input to both neural networks is the first-stage decision $\vect x$, the scenario embedding $\vect{g} \left( \bxi, \bpi; \Phi_{3}, \Phi_{4}\right)$ can be pre-computed for all $\bxi\in\Xi_O' \cup \Xi_F'$ before solving the surrogate main problem, hence there is no need to encode the MILP representations of the MLPs $\Phi_{3}$ and $\Phi_{4}$. {MLP $\Phi_1$ is encoded once per variable and  $\Phi_2$ is encoded once in total}. %
On the other hand, MLPs $\Phi_5$ and $\Phi_6$ are encoded multiple times, once for each scenario, so they must be small in size to limit the number of additional integer variables and constraints. Applying the notation for our network structure, the main problem can be reformulated as follows:        
\begin{subequations}
\label{eq:mp_neur2ro_full}
\begin{align}
    \min_{\substack{\vect x\in\gX,\vect y\in\gY, \Hat{f} \in \R, \\\bxi^{(a)}\in\Xi_O', \bxi^{(b)}\in\Xi_F'}} & 
        \quad \vect c (\bpi)^\intercal \vect x + \vect d(\bxi^{(a)},\bpi)^\intercal \vect y \\
    \text{s.t.} \qquad & 
        \quad W(\bxi^{(a)},\bpi) \vect y + T(\bxi^{(a)},\bpi) \vect x \le \vect h(\bxi^{(a)},\bpi), \\ 
        & \quad W(\bxi^{(b)},\bpi) \vect y + T(\bxi^{(b)},\bpi) \vect x \le \vect h(\bxi^{(b)},\bpi),  \\ 
        & \quad \Hat{f} = \vect{f}(\vect x, \bpi; \Phi_1, \Phi_2), \label{eq:mp_emb_1}\\
        & \quad \bxi^{(a)} \in {\arg\max}_{\bxi \in \Xi_O'} \left\{\Phi_5\left(\text{concat}\left[\Hat{f}, \Bar{\vect g}(\bxi, \bpi)\right]\right)  \right\}, \label{eq:mp_emb_2} \\
        & \quad \bxi^{(b)} \in {\arg\max}_{\bxi \in \Xi_F'} \left\{\Phi_6\left(\text{concat}\left[\Hat{f}, \Bar{\vect g}(\bxi, \bpi)\right]\right)\right\}, \label{eq:mp_emb_3}
\end{align}
\end{subequations}
where $\Bar{\vect g}(\bxi, \bpi) = \vect g(\bxi, \bpi; \Phi_3, \Phi_4)$ is pre-computed by evaluating the output of the corresponding networks via a forward pass. From constraint~\eqref{eq:mp_emb_1}, {it is clear that the MLPs that parameterize $\vect f$ are encoded once per variable for $\Phi_1$ and once in total for $\Phi_2$, respectively.} Constraints~\eqref{eq:mp_emb_2} and \eqref{eq:mp_emb_3} show that $\Phi_5$ and $\Phi_6$ are encoded for each scenario.

\subsubsection{Adversarial problems}
As for the surrogate adversarial problems~\eqref{eq:ap_o} and~\eqref{eq:ap_f}, the first-stage decision $\vect{x}^\star$ is fixed and so is its embedding $\vect{f} \left( \vect{x}^\star, \bpi; \Phi_{1}, \Phi_{2}\right)$. We pre-compute this embedding and maximize the output of the relevant neural network over the scenario space. Problem~\eqref{eq:ap_o} then writes as:
\begin{subequations}
\label{eq:neur2ro_ap_full}
\begin{align}
\max_{\bxi \in \Xi, \Hat{g} \in \R} & \quad \Phi_5\big(\text{concat}[\Bar{\vect f}(\vect x^\star, \bpi), \Hat{g}]\big), \\
\text{s.t.} & \quad \Hat{g} = \vect g (\bxi, \bpi; \Phi_3, \Phi_4),
\end{align}
\end{subequations} 
where $\Phi_3, \Phi_4$, and $\Phi_5$ are encoded only once each and $\Bar{\vect f}(\vect x^\star, \bpi) = \vect f(\vect x^\star, \bpi; \Phi_3, \Phi_4)$ is pre-computed by evaluating the output of the corresponding networks via a forward pass. Problem~\eqref{eq:ap_f} is reformulated similarly with $\Phi_6$ instead of $\Phi_5$.

\subsubsection{Encoding the $\text{argmax}$}
\label{sec:bigm}
The $\argmax$ operation over a set of neural network predictions can be implemented using the following standard modeling trick\footnote{\url{https://docs.mosek.com/modeling-cookbook/mio.html}, see `9.1.8 Maximum'}, applied for example to~\eqref{eq:mp_argmax1_neur2ro_} after introducing the variables $z_{\bxi} \in\{0,1\}, {\forall\bxi \in \Xi_O'}$:
\begin{equation}
\NN^O_{\Theta^{\star}}(\bpi,\vect x, \bxi) \leq t \leq \NN^O_{\Theta^{\star}}(\bpi,\vect x, \bxi) + M(1-z_{\bxi}), \; {\forall\bxi \in \Xi_O'}, \;
\sum_{\bxi \in \Xi_O'} z_{\bxi} = 1, \;
\bxi^{(a)} = \sum_{\bxi \in \Xi_O'} z_{\bxi}\bxi.
\end{equation}
Here, $M$ is a sufficiently large constant.
One of the scenarios that achieve the maximum neural network output across all scenarios in the set $\Xi_O'$ would see its indicator variable $z_{\bxi}$ set to one and thus $\bxi^{(a)} = \bxi$. The same encoding applies to the other $\argmax$ in~\eqref{eq:mp_argmax2_neur2ro_}.

\subsubsection{Encoding size} 
How many additional variables and linear constraints does one need in~\eqref{eq:mp_neur2ro_full} to encode all the relevant MLPs in the main problem? As our implementation uses the classical big-$M$ formulation for ReLU MLPs, we will consider the size of that formulation. For an MLP with $D$ hidden layers and $P$ neurons per layer, the number of binary variables and linear constraints grows as $\mathcal{O}(DP)$ (see Section 4.2.1 of~\citet{huchette2023deep}). For the surrogate master problem, this translates into $\mathcal{O}(nD_1P_1 + D_2P_2 + |\Xi_O'|D_5P_5 + |\Xi_F'|D_6P_6)$ where $D_k$ and $P_k$ are the depth and width of an MLP $\Phi_k$, respectively. In practice, we choose these MLP hyperparameters to be small and constant across instance sizes of the same problem. 
{In facility location, for example, $P_1=16, P_2=32, P_5=8, P_6=8$, and all MLPs have a single hidden layer, i.e., $D_k=1, k=1,\dots,6$.}
Treating these hyperparameters as constant, the surrogate main problem at a given iteration of~\Cref{alg:CCG_neur2ro} on a facility location instance has 
{$\mathcal{O}(16n + 8|\Xi_O'| + 8|\Xi_F'|)$} 
binary variables and linear inequalities in addition to the first-stage decision variables.
{In comparison, the number of integer variables of the original main problem~\eqref{eq:mp_levelset} grows as $\mathcal{O}(nm|\Xi'|)$ and the number of linear inequalities with $\mathcal{O}((n+m)|\Xi'|)$, where $n$ is the number of facilities and $m$ is the number of customers.} The surrogate main problem does not suffer from this bilinear growth in the number of second-stage variables \textit{and} scenarios and is only~\textit{linear} in these parameters. This keeps it tractable even as problem size and iterations grow.

As for the adversarial problems~\eqref{eq:neur2ro_ap_full}, they require encoding $\Phi_3, \Phi_4$, and $\Phi_5$ or $\Phi_6$ depending on whether one is trying to expand $\Xi_O'$ or $\Xi_F'$. Again, because all of the MLPs we use are small in size, this makes for relatively small MILPs. Compared to the vanilla CCG where a bilevel problem has to be solved, this formulation is much smaller.

It is worth noting that the big-$M$ formulation of ReLU MLPs has a weak linear relaxation. Stronger extended formulations have been proposed that could be explored in the context of the~\method{} method in the future~\citep{huchette2023deep}.

\subsubsection{Heuristic methods for the adversarial problem}
\label{sec:sub_grads}
The adversarial problems are MILPs when the uncertainty set can be modeled by linear constraints, e.g., if we consider polyhedral uncertainty sets.  The natural choice is to use an off-the-shelf MILP solver for optimizing the corresponding adversarial problems. Alternatively, one can use a faster heuristic to accelerate the learning-augmented CCG. Problems~\eqref{eq:ap_o} and \eqref{eq:ap_f} are maximizing the output of a neural network over the uncertainty set $\Xi$. A very similar task in the deep learning literature is that of ``adversarial attacks'' on trained neural networks~\citep{kurakin2017adversarial,madry2018towards}, where one seeks to maximize the output of an incorrect class prediction by perturbing an input feature vector in a certain convex norm-ball.

As a function of its inputs, a ReLU neural network has well-defined subgradients. Subgradients can be computed automatically using the auto-differentiation feature of modern deep-learning packages. An off-the-shelf non-linear programming solver such as Ipopt~\citep{wachter2006implementation} can be used as a heuristic for Problems~\eqref{eq:ap_o} and \eqref{eq:ap_f}. Ipopt uses the subgradients of the ReLU network as first-order information; second-order information is not used as ReLU networks have zero Hessian almost everywhere. We will compare an off-the-shelf MILP solver to Ipopt for solving the adversarial problem.  Note that other heuristics for optimizing over a trained neural network, based on sampling \citep{perakis2022optimizing} or local search \citep{tong2024optimization}, have been proposed.

\subsection{Ensuring (approximate) feasibility}
\label{sec:no_good_cuts}
Generally, \method{} may terminate with an infeasible solution.  To mitigate infeasibility, we can utilize standard techniques from mathematical optimization.  In the case of binary/integer first-stage variables, no-good cuts \citep{fischetti2016intersection} can be used to remove a candidate solution from the feasible region if an infeasible first-stage decision, $\vect x^\star$, is returned by \method{}.  
In general, this is intractable as it requires solving the adversarial problem; however, it can be efficiently approximated by checking if the second-stage problems are feasible for a set of $N$ scenarios sampled from the uncertainty set, as detailed in Algorithm~\ref{alg:feas_ngc}.

\begin{algorithm}[htbp!]
\footnotesize
\nonl\textbf{Input:} 
 $\bpi$, a problem instance's certain parameters;
 $\vect x^\star$, a solution;
 $\Xi$, the uncertainty set;
 $N$, a number of scenario samples.\\
\nonl\textbf{Output}: \textbf{True} if $\vect x^\star$ is feasible for all $N$ scenarios, \textbf{False} otherwise.
\begin{algorithmic}[1]
\STATE $\Xi_\text{feas}=\left\{ \bxi^{(j)} \overset{U}{\sim}\Xi\right\}_{j=1}^{N}$, a set of $N$ scenarios drawn from $\Xi$
    \FORALL{ $\left\{ \bxi^{(j)} \right\}_{j=1}^{N}$}
        \STATE Solve the second-stage problem:
        $\min_{\vect y\in\gY}\{\vect d(\bxi^{(j)},\bpi)^\intercal \vect y : W(\bxi^{(j)},\bpi) \vect y \le \vect h(\bxi^{(j)},\bpi) - T(\bxi^{(j)},\bpi )\vect x^\star \}$
        \IF{second-stage problems is infeasible}
        \STATE Add the no-good cut ($ \sum_{i:x_i^\star=0} x_i + \sum_{i:x_i^\star=1} (1-x_i) \geq 1 $) to the main problem
        \STATE return \textbf{False}
        \ENDIF
    \ENDFOR
    \STATE return \textbf{True}
 \caption{No-good Cuts Feasibility Improvement}
 \label{alg:feas_ngc}
 \end{algorithmic}
\end{algorithm}

Algorithm~\ref{alg:feas_ngc} can be seamlessly integrated within Algorithm~\ref{alg:CCG_neur2ro} as an additional termination criterion. If no cut is added, Algorithm~\ref{alg:CCG_neur2ro} terminates as usual. Otherwise, it continues.  This approach aims to ensure more feasible solutions are found at a relatively low additional cost in computing time.  Moreover, for problems where checking feasibility is tractable, the sampling in Algorithm~\ref{alg:feas_ngc} can even be replaced by exact approaches, ensuring \method{} computes a feasible solution.

\section{Theoretical Guarantees}
\label{sec:theory}

\subsection{Finite convergence}
Since our algorithm does not exactly apply the steps of the standard CCG, the convergence guarantee from the classical algorithm of~\citet{zeng2013solving} does not hold. Fortunately, we prove that finite convergence of the learning-augmented CCG~\PCref{alg:CCG_neur2ro} holds if there are finitely many first-stage solutions. This is the case if all first-stage variables are integer and $\gX$ is bounded, a property that holds for all applications considered in this work.
\begin{theorem}\label{thm:convergence}
If $\gX$ is finite, then~\Cref{alg:CCG_neur2ro} terminates after a finite number of iterations.
\end{theorem}
\begin{proof}
    For a pre-defined feasibility tolerance $\varepsilon>0$~\Cref{alg:CCG_neur2ro} moves on to the next iteration if either one of the following two conditions for adding a violating scenario is satisfied (lines 7 and 11):
\begin{equation}\label{eq:check_optimality}
\NN^O_{\Theta^{\star}}(\bpi, \vect x^\star, \bxi_O) \ge \max_{\bxi\in \Xi_O'} \NN^O_{\Theta^{\star}}(\bpi,\vect x^\star, \bxi) + \varepsilon,
\end{equation}
\begin{equation}\label{eq:check_feasibility}
\NN^F_{\Theta^{\star}}(\bpi, \vect x^\star, \bxi_F) \ge \varepsilon.
\end{equation}
We have to show that Conditions \eqref{eq:check_optimality} and \eqref{eq:check_feasibility} can only be true in a finite number of iterations. Since we terminate the algorithm if both conditions are false, the finite termination of the algorithm follows. If the solver for the adversarial problem \eqref{eq:ap_o} is consistent in that it always returns the same worst-case scenario for a given first-stage decision $\vect{x}^\star$, call it $\bxi_O(\vect x^\star)$, then the result is trivial: after adding this scenario to $\Xi_O'$, the next time the same solution $\vect x^\star$ is optimal for the master problem, the solution to the adversarial problem~\eqref{eq:ap_o} will be the exact same scenario $\bxi_O(\vect x^\star)$. Condition \eqref{eq:check_optimality} will not be satisfied. The same holds for the feasibility condition \eqref{eq:check_feasibility}. Note that the consistency of the solver can be guaranteed by maintaining a set of iterates $\vect{x}^\star$ along with their corresponding worst-case scenarios (i.e., the solutions to~\eqref{eq:check_optimality} and~\eqref{eq:check_feasibility}). Given a new first-stage solution, we check against this set and return its worst-case scenario if it is in the set.
\end{proof}

\subsection{Approximation guarantee}
Next, we analyze the approximation guarantee that our algorithm achieves. Since the (objective) neural network predictions are usually erroneous, we cannot guarantee an exact optimal solution. However, we can bound the optimality gap of a solution returned by~\method{} under the assumptions below. Denote the optimal second-stage value for a given instance $\bpi$, first-stage solution $\vect x$ and scenario $\bxi$ as
$$
\text{val} ( \vect x, \bxi, \bpi):= \min_{\vect y\in\gY}\{\vect d(\bxi,\bpi)^\intercal \vect y : W(\bxi ,\bpi) \vect y \le \vect h(\bxi,\bpi) - T(\bxi,\bpi )\vect x^\star \}
$$

We assume we have a given instance $\bpi\in \Pi$ to be solved. The assumptions are:
\begin{enumerate}[label=(\alph*)]
    \item for the given instance $\bpi$, the absolute error of the (objective) neural network is bounded by $\alpha >0$, i.e., 
        \begin{equation}\label{eq:prediction_error}
         | \NN^O_{\Theta^{\star}}(\vect x, \bxi, \bpi) - \text{val} ( \vect x, \bxi, \bpi )  |\le \alpha \quad \forall \vect x\in \gX, \bxi\in \Xi;
        \end{equation}
    \item relatively complete recourse holds, i.e., for every first-stage decision $\vect x\in \mathcal{X}$ and every scenario $\bxi\in \Xi$, there exists a feasible second-stage solution.
\end{enumerate}

\begin{theorem}\label{thm:approx_guarantee}
Let $\text{opt}$ be the optimal value of \eqref{eq:2s_ro} for a given instance $\bpi$. Under assumptions (a--b),~\method{} finds a solution $\vect x_{\NN}\in \gX$ with objective value
\[
\vect c^\intercal \vect x_{\NN} + 
\max_{\bxi\in\Xi} \; \text{val} ( \vect x_{\NN}, \bxi, \bpi ) \le \text{opt} + 2\alpha + \varepsilon.
\]
\end{theorem}
\begin{proof}
Note that since we assume relatively complete recourse, there is no need for the feasibility prediction network $\NN^F$ and the corresponding constraints \eqref{eq:mp_consfeas_neur2ro_} and~\eqref{eq:mp_argmax2_neur2ro_} can be removed from the surrogate main problem. Furthermore, in the following, we drop the parameter $\bpi$ from all expressions for ease of notation. 

Let $\Xi_O'$ be the set of scenarios generated during the execution of~\Cref{alg:CCG_neur2ro}. We denote by $\bxi^{(a)}(\vect x_{\NN})$ the argmax scenario selected in \eqref{eq:mp_argmax1_neur2ro_} for solution $ \vect x_{\NN}$ in the last iteration of the algorithm, i.e., 
\[
\bxi^{(a)}(\vect x_{\NN}) \in \argmax_{\bxi \in \Xi_O'} \;\NN^O_{\Theta^{\star}} ( \vect x_{\NN}, \bxi ).
\]
Furthermore, denote the optimal solution to \eqref{eq:2s_ro} as $\vect x_{\text{OPT}}$ and let $\bxi^{(a)}(\vect x_{\text{OPT}})$ be any argmax scenario selected in \eqref{eq:mp_argmax1_neur2ro_} for solution $\vect x_{\text{OPT}}$ in the last iteration of the algorithm. For every $\vect x\in \gX$, denote by $\bxi_{\text{OPT}}(\vect x)$ an arbitrary worst-case scenario of the exact two-stage problem, i.e.,
\begin{equation}\label{eq:def_worst-case_scenario}
\bxi_{\text{OPT}}(\vect x) \in \argmax_{\xi\in \Xi} \;\text{val} ( \vect x, \bxi ).
\end{equation}

Then the following inequalities hold:
\begin{align*}
    & \vect c^\intercal \vect x_{\NN} + \max_{\bxi\in\Xi} \;\text{val} ( \vect x_{\NN}, \bxi ) - \text{opt} \\
    & \qquad = \vect c^\intercal \vect x_{\NN} + \text{val} ( \vect x_{\NN}, \bxi_{\text{OPT}} (\vect x_{\NN}) ) - \vect c^\intercal \vect x_{\text{OPT}} - \text{val} ( \vect x_{\text{OPT}}, \bxi_{\text{OPT}}(\vect x_{\text{OPT}}) ) \\
    & \qquad \le \vect c^\intercal \vect x_{\NN} +  \NN_{\Theta^{\star}}^O ( \vect x_{\NN}, \bxi_{\text{OPT}} (\vect x_{\NN}) ) - \vect c^\intercal \vect x_{\text{OPT}} - \text{val} ( \vect x_{\text{OPT}}, \bxi_{\text{OPT}}(\vect x_{\text{OPT}}) ) + \alpha && \text{by~\eqref{eq:prediction_error}}\\
    & \qquad \le \vect c^\intercal \vect x_{\NN} +  \NN_{\Theta^{\star}}^O ( \vect x_{\NN}, \bxi_{\text{OPT}} (\vect x_{\NN}) ) - \vect c^\intercal \vect x_{\text{OPT}} - \text{val} ( \vect x_{\text{OPT}}, \bxi^{(a)}(\vect x_{\text{OPT}}) ) + \alpha && \text{by~\eqref{eq:def_worst-case_scenario}}\\
    & \qquad \le \vect c^\intercal \vect x_{\NN} +  \NN_{\Theta^{\star}}^O ( \vect x_{\NN}, \bxi^{(a)} (\vect x_{\NN}) ) - \vect c^\intercal \vect x_{\text{OPT}} - \text{val} ( \vect x_{\text{OPT}}, \bxi^{(a)}(\vect x_{\text{OPT}}) ) + \alpha + \varepsilon \\
    & \qquad \le \vect c^\intercal \vect x_{\NN} +  \NN_{\Theta^{\star}}^O ( \vect x_{\NN}, \bxi^{(a)} (\vect x_{\NN}) ) - \vect c^\intercal \vect x_{\NN} - \text{val} ( \vect x_{\NN}, \bxi^{(a)}(\vect x_{\NN}) ) + \alpha + \varepsilon \\
    & \qquad \le  \NN_{\Theta^{\star}}^O(\vect x_{\NN}, \bxi^{(a)} (\vect x_{\NN}) ) - \NN_{\Theta^{\star}}^O ( \vect x_{\NN}, \bxi^{(a)}(\vect x_{\NN}) ) + 2\alpha + \varepsilon && \text{by~\eqref{eq:prediction_error}}\\
    & \qquad = 2\alpha + \varepsilon,
\end{align*}
where the third inequality holds because condition \eqref{eq:check_optimality} is not fulfilled in the last iteration of the algorithm and the fourth inequality holds since $\vect x_{\NN}$ is the optimal solution of the surrogate main problem in the last iteration.
\end{proof}

We now consider the case where relatively complete recourse does not hold, i.e., it may be possible that for a solution $\vect x\in\gX$ there exists a scenario $\bxi \in \Xi$ such that the second-stage problem is infeasible. In this case, we use the feasibility prediction network $\NN^F$. When using $\NN^F$,~\Cref{alg:CCG_neur2ro} may return solutions with an arbitrarily large optimality gap. This is because if $\NN^F$ misclassifies any good solution as infeasible for at least one scenario, then these solutions are cut off by constraint \eqref{eq:mp_consfeas_neur2ro_} whereas the feasible solutions can be suboptimal. It could even happen that our algorithm predicts every solution to be infeasible and does not return any solution, although the problem is feasible. Hence, in the most general case, no approximation guarantee similar to the one in Theorem \ref{thm:approx_guarantee} can be derived.

\section{Experiments}
\label{sec:experiments}

All experiments were run on a computing cluster with an Intel Xeon CPU E5-2683 and Nvidia Tesla P100 GPU with 64GB of RAM (for training).  Pytorch 1.12.1 \citep{NEURIPS2019_9015} was used for all learning models.  Gurobi 10.0.2 \citep{gurobi} was used as the MILP solver and gurobi-machinelearning 1.3.0 was used to embed the neural networks into MILPs.  For evaluation, all solving was limited to 3 hours. Our code is available at \url{https://github.com/khalil-research/Neur2RO}.

\subsection{Benchmark Problems}
In the following sections, we evaluate \method{} on three two-stage robust optimization problems. \Cref{app:benchmarks} provides complete formulations and definitions of uncertainty sets, as well as details of instance generation. Here, we briefly outline the benchmark problems and relevant parameters used in the numerical experiments.

\begin{itemize}
    \item \textbf{Two-stage robust knapsack problem} proposed by \citet{arslan2022min} and inspired by \citet{ben2009robust}.  We evaluate the publicly available test instances\footnote{Instances and results available at \url{https://github.com/borisdetienne/RobustDecomposition}} from \citet{arslan2022min}. The instances are categorized into four groups depending on the correlation of the nominal profits of the items with their costs:  uncorrelated (UN), weakly correlated (WC), almost strongly correlated (ASC), and strongly correlated (SC); more correlated instances are much harder to solve. We consider instances with $n\in\{20,30,40,50,60,70,80\}$ items. These instances contain only objective uncertainty.
    
    \item \textbf{Capital budgeting problem} proposed by \citet{hanasusanto2015k}. These problem instances are generated using the procedure presented by \citet{subramanyam2020k}.  We consider instances with $n\in\{10,20,30,40,50\}$ items containing objective and constraint uncertainty.
    
    \item \textbf{Two-stage robust facility location problem} with binary second-stage variables as defined in \eqref{eq:FL2}. Compared to the formulation described in the introduction, we have demand uncertainty (constraint uncertainty) as defined in \eqref{eq:FL} and additionally, disruption uncertainty (objective uncertainty) as defined in \eqref{eq:FL2}. These instances require feasibility prediction there is not relatively complete recourse.   
    We evaluate instances with $n\in\{5,10,20\}$ facilities and $m\in\{10,20,50\}$ customers, for $m \geq n$.
\end{itemize}

\subsection{Baseline methods}
\begin{itemize}
    \item \textbf{Static robust optimization} (\static{}): A simplification of 2RO based on a single-level reformulation of the problem, where all second-stage decisions are moved to the first level such that no adaptivity to the uncertainty is possible. In \citet{ben2009robust}, a computationally tractable formulation is given to deal with such problems. This is equivalent to but more efficient than $k$-adaptability with $k=1$. Details of the \static{} models for each benchmark are provided in~\Cref{app:static}
    \item \textbf{Branch-and-price} (BP): Proposed by \citet{arslan2022decomposition}, BP is an exact algorithm for 2RO problems with only objective uncertainty, e.g., the knapsack problem.  
    \item $k$\textbf{-adaptability}: Another simplification of 2RO, where $k$ second-stage reactions are already calculated in the first-stage. This approach provides heuristic solutions to the 2RO problem, where larger values of $k$ increase the quality of the solution. Since the capital budgeting and facility location problems have constraint uncertainty, we use the $k$-adaptability branch-and-bound algorithm of \cite{subramanyam2020k}, with $k=2,5,10$ as a baseline. We do not include $k$-adaptability for the knapsack problem as \cite{arslan2022decomposition} demonstrate that BP significantly outperforms $k$-adaptability.
\end{itemize}

\subsection{Learning Algorithms} 
We consider three variants of our approach. 
The first solves the adversarial problem in \method{} as a MILP (\methodmip{}).
The second solves the adversarial problem using Ipopt (\methodipo{}) as discussed in Section~\ref{sec:sub_grads}. 
The third solves the adversarial as a MILP and utilizes Algorithm~\ref{sec:no_good_cuts} with $N=100$ sampled scenarios 
(\methodngc{}).
The latter is only applicable to problems with constraint uncertainty.  
For all variants, we terminate a solve of the main or adversarial problems early if no improvement in the solution is observed in 180 seconds.  For Ipopt, we limit the number of iterations to 50 per adversarial solve. 

Details regarding model hyperparameters, features, and times associated with data collection and training are reported in Appendix~\ref{app:ml_details}. One model is trained for each benchmark problem, where in the worst case, this takes~$\sim4.5$ hours, which is relatively insignificant considering baseline algorithms often take more than 3 hours to terminate on a single instance.

\subsection{Evaluation}
After each method returns a first-stage decision, we obtain the corresponding objective value by solving \eqref{eq:2s_ro} for a fixed first-stage decision.  However, this is handled separately for problems with objective and constraints uncertainty. For problems with only objective uncertainty, this is done exactly using constraint generation, while for problems with constraint uncertainty, we approximate the worst-case objective via sampling. Our experiments sample 10,000 scenarios and additionally use all worst-case scenarios computed via each algorithm. Details of each procedure are provided in Appendix~\ref{app:ove}.

For evaluation of the efficacy and efficiency of \method{}, we report the following metrics.  
\begin{itemize}
    \item \textbf{Relative Error}:  We use the relative error to evaluate the solution quality, i.e., the objective gap to the best-known solution achieved by any algorithm on an instance.  The relative error of an algorithm $\mathcal{A}$ is computed as $100\cdot\frac{| obj^\star - obj(\mathcal{A}) |}{|obj^\star|}$, where $obj^\star$ is the best known solution.
    
    \item \textbf{Solving Time}: The solving time is the time it takes for each algorithm to terminate.  Note that the limit is $\sim10,800$ seconds as a 3-hour time limit is given for each instance.  

    \item \textbf{Feasibility}: Only applicable to the facility location problem, we evaluate the fraction of instances for which a feasible solution was found. For $k$-adaptability, whenever a solution is returned, it is guaranteed to be feasible.  
    For all variants of \method{}, a computed solution may be infeasible. 
    As such, an instance is reported as feasible if the solution has feasible second-stage problems for all of the sampled scenarios. This is not reported for capital budgeting as the baselines are always feasible, and ensuring feasibility with \method{} is trivial, as discussed in~\Cref{app:cb_trivial}.
\end{itemize}   
In addition, to provide a more comprehensive evaluation, we compare the solution quality of each baseline at the termination time of \method{}. These results are reported in~\Cref{app:ml_times}.  

\subsection{Numerical Results}
\begin{figure}
    \centering
    \begin{minipage}{0.5\textwidth}
        \centering
        \includegraphics[width=0.99\textwidth]{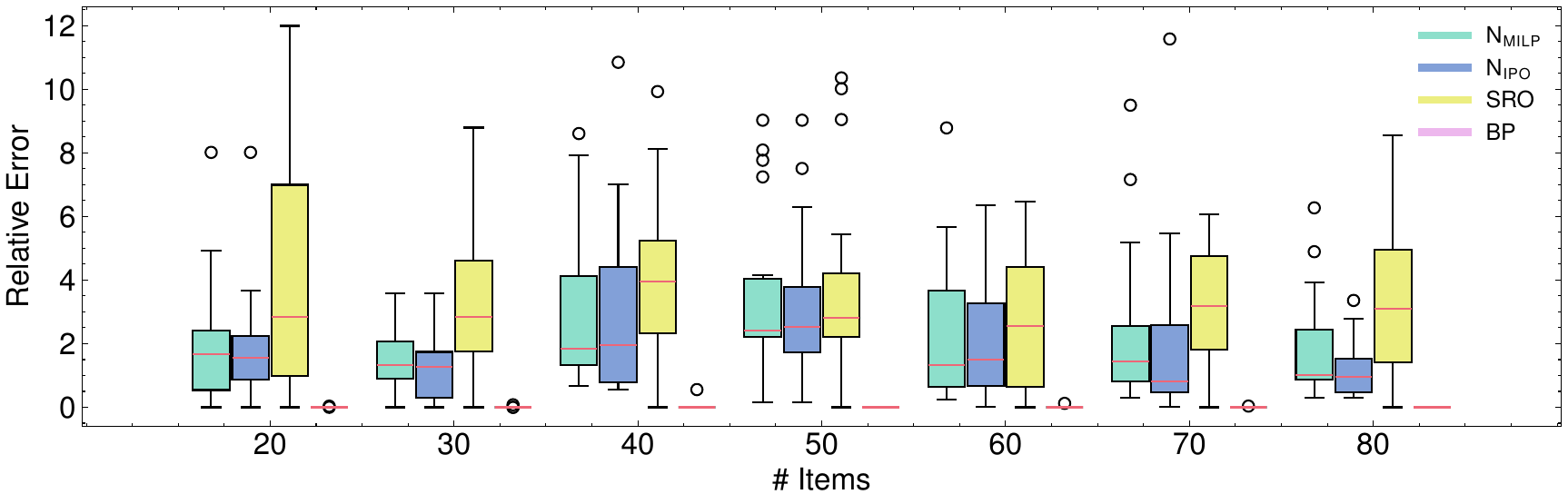} 
    \end{minipage}\hfill 
    \begin{minipage}{0.5\textwidth}
        \centering
        \includegraphics[width=0.99\textwidth]{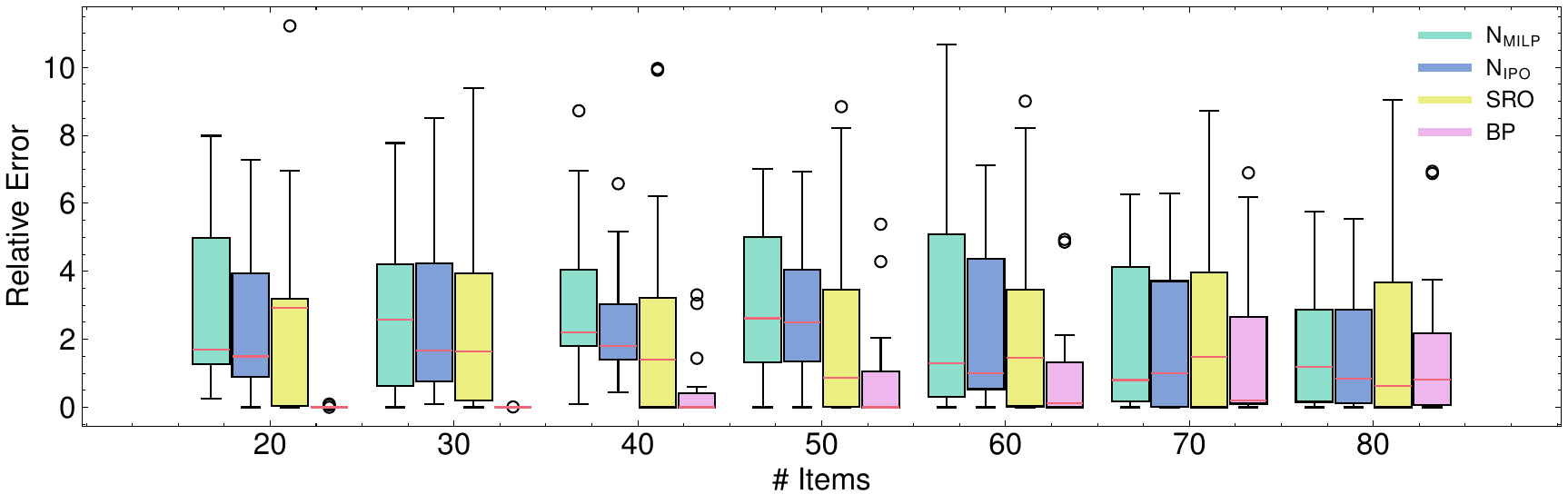} 
    \end{minipage}\hfill \\
    \centering
    \begin{minipage}{0.5\textwidth}
        \centering
        \includegraphics[width=0.99\textwidth]{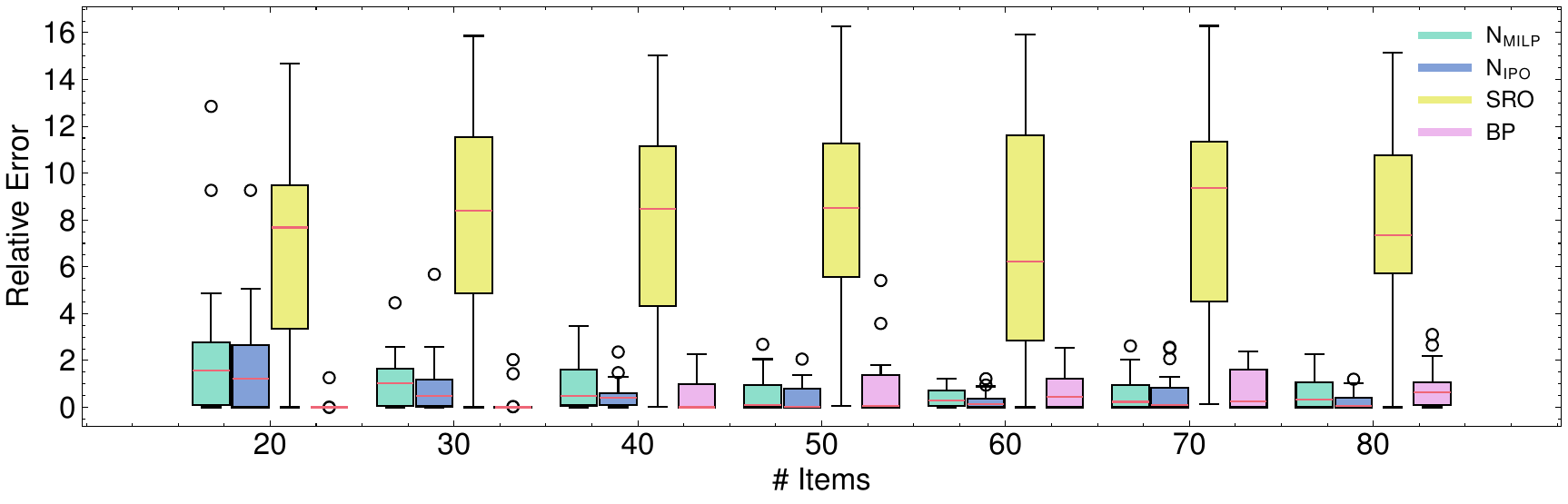} 
    \end{minipage}\hfill 
    \begin{minipage}{0.5\textwidth}
        \centering
        \includegraphics[width=0.99\textwidth]{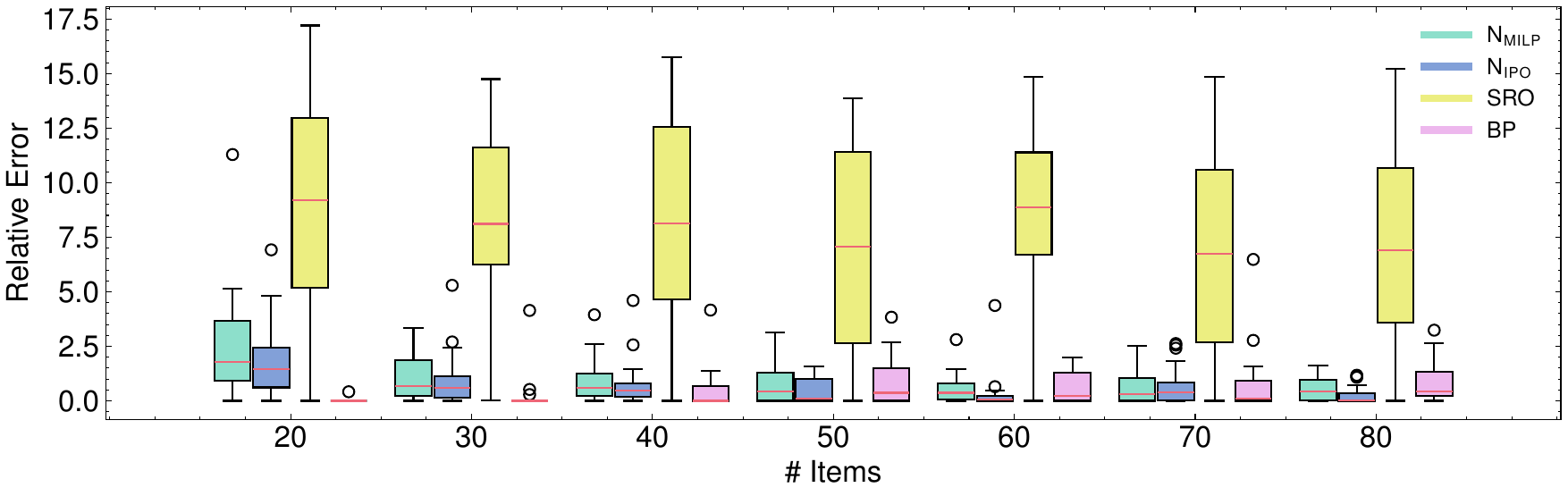} 
    \end{minipage}\hfill
    \caption{Box plot of relative errors for \methodmip{}, \methodipo{}, \static{}, and BP on knapsack instances.  Each box contains the relative errors over 18 instances.  The uncorrelated, weakly correlated, almost strongly correlated, and strongly correlated instances are in the top-left, top-right, bottom-left, and bottom-right, respectively.}
    \label{fig:kp_boxplot}
\end{figure}

\begin{figure}
    \centering
    \begin{minipage}{0.6\textwidth}
        \centering
        \includegraphics[width=0.9\textwidth]{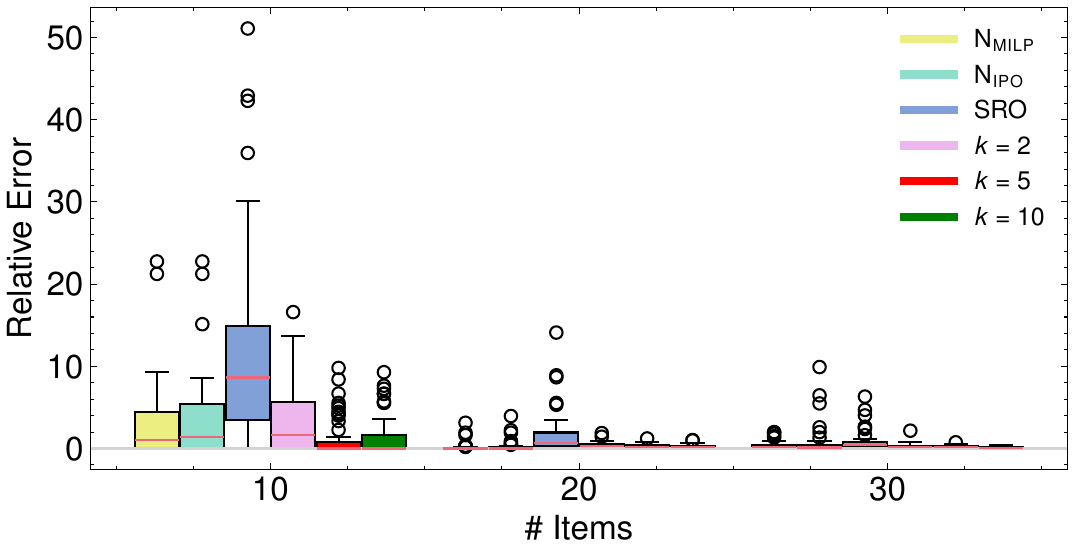} 
    \end{minipage}\hfill
    \begin{minipage}{0.4\textwidth}
        \centering
        \includegraphics[width=0.9\textwidth]{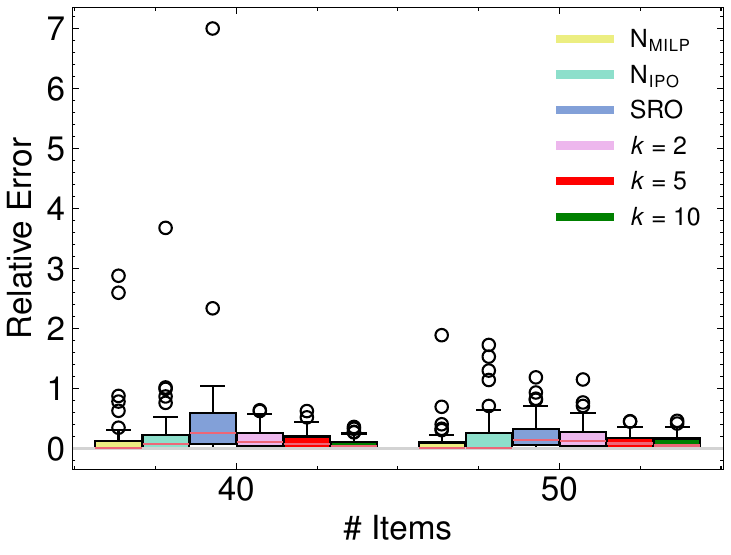} 
    \end{minipage}\hfill
    \caption{Box plot of relative errors for \methodmip{}, \methodipo{}, \static{}, 
    and $k$-adaptability on capital budgeting instances.  Each box contains the relative errors over 50 instances.}
    \label{fig:cb_boxplot}
\end{figure}

\begin{figure}[htbp!]
\centering
\includegraphics[width=1.0\textwidth]{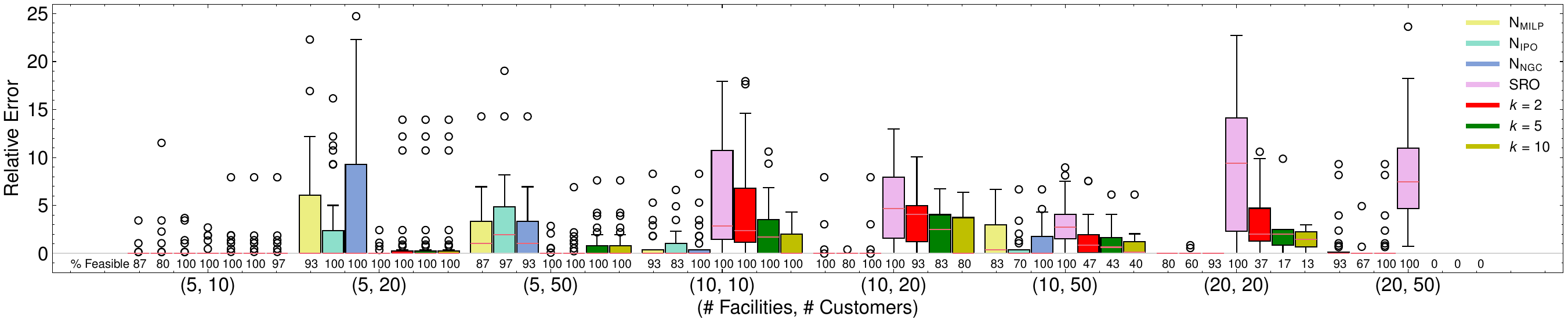}
\caption{
Box plot of relative errors for \methodmip{}, \methodipo{}, \methodngc{}, \static{}, and $k$-adaptability on facility location instances.  Only feasible solutions are used to compute the relative error, and the percentage of possible solutions for each instance size and method are provided below the respective boxplot.  For algorithms that do not compute feasible solutions for any of the instances, no boxes are plotted.  Each box contains the relative errors over 30 instances excluding infeasible cases.
}
\label{fig:cflp_boxplot}
\end{figure}

In Tables~\ref{tab:kp_results}--\ref{tab:cflp_results}, the median relative errors and solving times are reported. Figures~\ref{fig:kp_boxplot}--\ref{fig:cflp_boxplot} provide a detailed view of the distribution of relative errors across instances via box plots. In the remainder of this section, we discuss the numerical results for each problem. 

\begin{table}[htbp!]\centering\resizebox{1.0\textwidth}{!}{
    \centering
    \small
    \begin{tabular}{l|c|rrrr|rrrrcl|c|rrrr|rrrr}
        \cmidrule[0.8pt]{0-9}
        \cmidrule[0.8pt]{12-21}
        Correlation  & \# items & \multicolumn{4}{c|}{Median Relative Error}             &  \multicolumn{4}{c}{Times}                  & \hspace{0.25cm}  &  Correlation  & \# items & \multicolumn{4}{c|}{Median Relative Error}             &  \multicolumn{4}{c}{Times}  \\
        Type         &          & \methodmip{}  & \methodipo{} & \static{} & BP   &  \methodmip{}  & \methodipo{} & \static{} & BP   & \hspace{0.25cm}  &  Type         &          & \methodmip{}  & \methodipo{} & \static{} & BP   & \methodmip{}  & \methodipo{} & \static{} & BP     \\               
        \cmidrule[0.6pt]{0-9}
        \cmidrule[0.6pt]{12-21}
        \multirow{7}{*}{Uncorrelated}                                                          
             & 20 & 1.683 & 1.564 & 2.850 & \textbf{0.000} & 3 & 1 & \textbf{0} & 0                             & & \multirow{7}{*}{\shortstack[l]{Almost\\ Strongly\\ Correlated}} & 20 & 1.556 & 1.229 & 7.678 & \textbf{0.000} & 4 & 1 & \textbf{0} & 9  \\
             & 30 & 1.335 & 1.278 & 2.838 & \textbf{0.000} & 4 & 1 & \textbf{0} & 1                             & &  & 30 & 1.041 & 0.474 & 8.398 & \textbf{0.000} & 6 & 1 & \textbf{0} & 2708  \\
             & 40 & 1.843 & 1.962 & 3.949 & \textbf{0.000} & 8 & 2 & \textbf{0} & 3                             & &  & 40 & 0.501 & 0.401 & 8.466 & \textbf{0.000} & 9 & 2 & \textbf{1} & 4744  \\
             & 50 & 2.426 & 2.539 & 2.824 & \textbf{0.000} & 7 & 7 & \textbf{0} & 12                            & &  & 50 & 0.108 & \textbf{0.030} & 8.503 & 0.076 & 8 & 2 & \textbf{2} & 8852  \\
             & 60 & 1.341 & 1.509 & 2.543 & \textbf{0.000} & 13 & 2 & \textbf{0} & 18                           & &  & 60 & 0.302 & \textbf{0.146} & 6.243 & 0.452 & 15 & 3 & \textbf{2} & 10261 \\
             & 70 & 1.433 & 0.815 & 3.192 & \textbf{0.000} & 16 & 4 & \textbf{0} & 46                           & &  & 70 & 0.232 & \textbf{0.105} & 9.349 & 0.239 & 17 & \textbf{2} & 4 & 10800  \\
             & 80 & 1.014 & 0.958 & 3.108 & \textbf{0.000} & 11 & 2 & \textbf{0} & 388                          & &  & 80 & 0.346 & \textbf{0.057} & 7.335 & 0.650 & 13 & \textbf{8} & 22 & 10800 \\
        \cmidrule[0.6pt]{0-9}
        \cmidrule[0.6pt]{12-21}
        \multirow{7}{*}{\shortstack[l]{Weakly\\ Correlated}}
            & 20 & 1.705 & 1.497 & 2.919 & \textbf{0.000} & 5 & 1 & \textbf{0} & 29                             & &  \multirow{7}{*}{\shortstack[l]{Strongly\\ Correlated}}   & 20 & 1.774 & 1.447 & 9.192 & \textbf{0.000} & 5 & 1 & \textbf{0} & 9 \\
            & 30 & 2.580 & 1.661 & 1.651 & \textbf{0.000} & 9 & 2 & \textbf{0} & 454                            & &    & 30 & 0.670 & 0.609 & 8.110 & \textbf{0.000} & 5 & 1 & \textbf{0} & 2473 \\
            & 40 & 2.193 & 1.808 & 1.405 & \textbf{0.000} & 19 & 5 & \textbf{0} & 6179                          & &    & 40 & 0.612 & 0.462 & 8.136 & \textbf{0.000} & 22 & 2 & \textbf{1} & 5665 \\
            & 50 & 2.617 & 2.503 & 0.865 & \textbf{0.000} & 32 & 2 & \textbf{0} & 8465                          & &    & 50 & 0.424 & \textbf{0.090} & 7.073 & 0.370 & 10 & 2 & \textbf{2} & 8240 \\
            & 60 & 1.294 & 1.008 & 1.451 & \textbf{0.116} & 78 & 10 & \textbf{0} & 9242                         & &    & 60 & 0.369 & \textbf{0.078} & 8.855 & 0.209 & 11 & 28 & \textbf{2} & 10800 \\
            & 70 & 0.802 & 0.995 & 1.488 & \textbf{0.194} & 16 & 50 & \textbf{0} & 10800                        & &    & 70 & 0.321 & 0.371 & 6.727 & \textbf{0.093} & 13 & \textbf{2} & 7 & 10800 \\
            & 80 & 1.180 & 0.831 & \textbf{0.634} & 0.812 & 33 & 74 & \textbf{0} & 10800                        & &    & 80 & 0.436 & \textbf{0.028} & 6.892 & 0.421 & 12 & 10 & \textbf{6} & 10800 \\
        \cmidrule[0.8pt]{0-9}
        \cmidrule[0.8pt]{12-21}
    \end{tabular}}
    \caption{Knapsack problem: median relative error and average solving times in seconds. For each row, the results are computed over 18 instances. The smallest (best) values in each row/metric are in bold.}
    \label{tab:kp_results}
\end{table}

\begin{table*}[htbp!]\centering\resizebox{0.65\textwidth}{!}{
\begin{tabular}{c|rrrrrr|rrrrrr}
\toprule
\# items & \multicolumn{6}{c|}{Median Relative Error} & \multicolumn{6}{c}{Times}  \\
 & \methodmip{} & \methodipo{} & \static{} & $k=2$ & $k=5$ & $k=10$ & \methodmip{}\ \  & \methodipo{}\ \  & \static{}\  & $k=2$\ \  & $k=5$\ \  & $k=10$\ \  \\
\midrule
10 & 1.014 & 1.383 & 7.945 & 1.579 & \textbf{0.000} & \textbf{0.000} & 60 & 2 & \textbf{0} & 20 & 9561 & 10800 \\
20 & \textbf{0.000} & \textbf{0.000} & 0.669 & 0.314 & 0.230 & 0.157 & 216 & 236 & \textbf{0} & 8702 & 10800 & 10800 \\
30 & 0.156 & 0.126 & 0.473 & 0.176 & 0.167 & \textbf{0.089} & 664 & 264 & \textbf{0} & 10801 & 10800 & 10800 \\
40 & \textbf{0.007} & 0.078 & 0.254 & 0.114 & 0.068 & 0.046 & 630 & 496 & \textbf{0} & 10806 & 10801 & 10801 \\
50 & \textbf{0.008} & 0.013 & 0.149 & 0.129 & 0.083 & 0.050 & 808 & 444 & \textbf{0} & 10807 & 10804 & 10801 \\
\bottomrule
\end{tabular}}
\caption{Capital budgeting problem: median relative error and average solving times.  For each row, the results are computed over 50 instances. The smallest (best) values in each row/metric are in bold.}
\label{tab:cb_results}
\end{table*}

\begin{table*}[htbp!]\centering\resizebox{1.0\textwidth}{!}{
\begin{tabular}{l|rrrrrrr|rrrrrrr|rrrrrrr}
\toprule
($m$, $n$) & \multicolumn{7}{c|}{Median RE} & \multicolumn{7}{c|}{Times (seconds)} & \multicolumn{7}{c}{Percent of feasible/found solution)} \\
 & \methodmip{} & \methodipo{} & \methodngc{} & \static & $k=2$ & $k=5$ & $k=10$ & \methodmip{}\ \  & \methodipo{}\ \  & \methodngc{}\ \  & \static\ \  & $k=2$\ \  & $k=5$\ \  & $k=10$\ \  & \methodmip{}\ \ \  & \methodipo{}\ \ \  & \methodngc{}\ \ \  & \static\ \ \  & $k=2$\ \ \  & $k=5$\ \ \  & $k=10$\ \ \  \\
\midrule
(5, 10) & \textbf{0.000} & \textbf{0.000} & \textbf{0.000} & \textbf{0.000} & \textbf{0.000} & \textbf{0.000} & \textbf{0.000} & 5 & 1 & 7 & \textbf{0} & 54 & 2739 & 3043 & 87 & 80 & \textbf{100} & \textbf{100} & \textbf{100} & \textbf{100} & \textbf{100} \\
(5, 20) & \textbf{0.000} & \textbf{0.000} & \textbf{0.000} & \textbf{0.000} & \textbf{0.000} & \textbf{0.000} & \textbf{0.000} & 12 & 1 & 14 & \textbf{0} & 174 & 2717 & 2906 & 93 & \textbf{100} & \textbf{100} & \textbf{100} & \textbf{100} & \textbf{100} & \textbf{100} \\
(5, 50) & 1.030 & 1.967 & 1.030 & \textbf{0.000} & \textbf{0.000} & \textbf{0.000} & \textbf{0.000} & 9 & 4 & 14 & \textbf{0} & 4653 & 8047 & 8281 & 87 & 97 & 93 & \textbf{100} & \textbf{100} & \textbf{100} & \textbf{100} \\
(10, 10) & \textbf{0.000} & \textbf{0.000} & \textbf{0.000} & 2.877 & 2.378 & 1.724 & \textbf{0.000} & 9 & 4 & 14 & \textbf{1} & 8139 & 10465 & 10465 & 93 & 83 & \textbf{100} & \textbf{100} & \textbf{100} & \textbf{100} & \textbf{100} \\
(10, 20) & \textbf{0.000} & \textbf{0.000} & \textbf{0.000} & 4.712 & 4.070 & 2.495 & \textbf{0.000} & 14 & 4 & 16 & \textbf{2} & 9722 & 10800 & 10800 & \textbf{100} & 80 & \textbf{100} & \textbf{100} & 97 & 87 & 80 \\
(10, 50) & 0.379 & \textbf{0.000} & \textbf{0.000} & 2.758 & 0.900 & 0.663 & 0.122 & 29 & \textbf{14} & 40 & 22 & 10548 & 10801 & 10800 & 83 & 70 & \textbf{100} & \textbf{100} & 47 & 43 & 40 \\
(20, 20) & \textbf{0.000} & \textbf{0.000} & \textbf{0.000} & 9.428 & 2.041 & 2.041 & 1.465 & 23 & \textbf{9} & 49 & 361 & 10800 & 10800 & 10800 & 80 & 60 & 93 & \textbf{100} & 37 & 17 & 13 \\
(20, 50) & \textbf{0.000} & \textbf{0.000} & \textbf{0.000} & 7.460 & - & - & - & \textbf{37} & 72 & 55 & 1344 & - & - & - & 93 & 67 & \textbf{100} & \textbf{100} & 0 & 0 & 0 \\
\bottomrule
\end{tabular}}
\caption{Facility location: median relative errors, average solving times in seconds, and feasibility rates. Symbol $m$ is the number of facilities and $n$ is the number of customers. For each row, the results are computed over 30 instances. The best values in each row/metric are in bold. A value of ``-'' indicates that no feasible solutions were found.}
\label{tab:cflp_results}
\end{table*}

\paragraph{Knapsack problem.} Table~\ref{tab:kp_results} and Figure~\ref{fig:kp_boxplot} demonstrate a clear improvement in scalability, with average solving times of \methodmip{} and \methodipo{} ranging from 1 to 78 seconds, while the solving times of BP increase with the difficulty of the instances, i.e., those with a large number of items and (almost) strong correlation. SRO is generally very efficient, but the solution quality is markedly worse than any other method. For the more difficult instances, \methodmip{} and \methodipo{} generally find better quality solutions and are over 100 times faster than BP, which is a very strong result considering BP is the state-of-the-art for problems with objective uncertainty.  Moreover, even for easy instances, both variants of \methodmip{} consistently compute much better solutions than SRO while taking significantly less time than BP.  

\paragraph{Capital budgeting.} Table~\ref{tab:cb_results} shows that \methodmip{} achieves the lowest median relative error for 20, 40, and 50-item instances, i.e., the two largest and most challenging instance sets.  Furthermore, \methodipo{} terminates in less time, while still achieving only marginally worse quality solutions. The distribution of the relative error for 40 and 50-item instances provided in Figure~\ref{fig:cb_boxplot} is consistent with the median result, as it illustrates that \methodmip{} finds high-quality solutions on the majority of the instances.  
In terms of solving time, \methodmip{} generally converges much faster than $k$-adaptability, resulting in a favorable trade-off: we can find better or equally good solutions 10 to 100 times faster.  Moreover, taking the incumbent solution found by $k$-adaptability at \methodmip{} termination time typically yields worse solutions than those reported in Table~\ref{tab:cb_results}; see~\Cref{app:ml_times} 
for details. SRO generally computes solutions much more efficiently than both other methods. However, the solution quality is much worse.

\paragraph{Facility location.} This is perhaps the most vital benchmark to demonstrate the efficacy of \method{} as it has the most complex constraints and two types of uncertainty. From~\Cref{tab:cflp_results} and~\Cref{fig:cflp_boxplot}, we can notice that $k$-adaptability fails to find solutions for larger instances. Moreover, all variants of \method{} find better solutions than SRO.  

\methodipo{} generally terminates more quickly than \methodmip{}, and achieves similar solution quality.  However, it results in more infeasible solutions.  
In terms of feasibility, \methodmip{} computes feasible solutions on $89\%$ of the instances, which justifies the feasibility prediction element of our approach.
Moreover, \methodngc{} further improves upon \methodmip{} by computing feasible solutions on  $98\%$ of the instances, with only marginal increases in runtime, demonstrating the efficacy of Algorithm~\ref{alg:feas_ngc},
and more generally, how solution refinement techniques can complement \method{}.  

In summary, for all three benchmark problems, \method{} achieves high-quality solutions.  For relatively easy or small instances, state-of-the-art methods sometimes find slightly better solutions, often at a much higher computational cost. However, as the instances become more difficult, \method{} demonstrates a clear improvement in overall solution quality and computing time. 

\section{Conclusion}
\label{sec:conclusion}
\method{} is the first learning-augmented approach for two-stage robust optimization with constraint uncertainty and integer decision variables. Key to the success of our method is the careful combination of a tailor-made neural network architecture with column-and-constraint generation. This intricate hybrid admits some theoretical convergence and approximation guarantees while performing significantly better than existing 2RO algorithms on widely used benchmark problems. Our work suggests several directions for future research. While empirically accurate, our neural network architecture is not guaranteed to be capable of representing the second-stage value function; exploring this representation question could be of theoretical interest. Another theoretical question relates to the sample complexity of our learning task, i.e., how many instances, first-stage decisions, scenarios, and the corresponding objective value must one collect to guarantee generalization to unseen instances with high probability? Similar questions in the learning-to-optimize space have been recently considered~\citep{balcan2021much}. {Generalizing from a two-stage to a }multi-stage adjustable robust optimization setting~\citep{lorca2016multistage,lappas2016multi} would also be a natural target for extending~\method.

\bibliography{main.bib}
\label{sec:bib}

\clearpage
\appendix 
\section{Benchmarks}
\label{app:benchmarks}
This section includes the MILP problem formulations and parameter generation/selection of the experimental results of the benchmarks and their clarifications.

\subsection{Knapsack Problem}
We consider the two-stage knapsack problem as defined in \cite{arslan2022decomposition} with a set of $n$ items. Each item $i$ has a weight $c_i$ and an uncertain profit $p_i(\bxi) = \Bar{p}_i - \xi_i \Hat{p}_i$, where $\Bar{p}_i$ is the expected profit, $\Hat{p}_i$ its maximum deviation and $\xi_i$ the \emph{uncertain} profit degradation factor, where the degradation happens after the first stage. In this problem we have a budgeted uncertainty set $\Xi = \{\bxi \in [0, 1]^n: \sum_{i=1}^n \xi_i \leq \Gamma\}$. The first-stage decision is to choose a subset of items to produce. Then in the second stage, there are three different responses to the profit degradation: (i) accept the degraded profit, (ii) repair the item by using an additional $t_i$ units from the budget to recover the original profit $\Bar{p}_i$, or (iii) outsource the item for a cost of $f_i$ units, such that the item's profit results in $\Bar{p}_i - f_i$. This gives the following problem formulation: 
\begin{subequations}  \label{KP}%
\begin{align}
    \min_{\vect x \in \{0, 1\}^n} \; \max_{\bxi \in \Xi} \; \min_{\vect y \in \{0, 1\}^n, \vect r \in \{0, 1\}^n} \quad & \sum_{i=1}^n (f_i - \Bar{p}_i)x_i + (\Hat{p}_i\xi_i - f_i)y_i - \Hat{p}_i \xi_i r_i \\
    \text{s.t.} \qquad & \sum_{i=1}^n c_i y_i + t_i r_i \leq C \\ 
    & r_i \le y_i \leq x_i \hspace{5cm} \forall i \in \{1, \ldots, n\},
\end{align}
\end{subequations}  
where $x_i$ is the first-stage decision to produce item $i$. For the second-stage decisions, we have $y_i$ and $r_i$: (i) $y_i = 1$ if item $i$ is produced \emph{without} repairing and $y_i = 0$ if the item is outsourced, and (ii) $r_i$ is the decision for repairing item $i$.

\paragraph{Instances in Experiments:} For the knapsack instances, we directly use the entire set of instances from \citet{arslan2022decomposition}, which are publicly available at \url{https://github.com/borisdetienne/RobustDecomposition}. \cite{arslan2022decomposition} describe full details of the instance generation procedure.

\subsection{Capital Budgeting Problem}
For the capital budgeting problem, we consider the formulation proposed by \citet{hanasusanto2015k} and \citet{subramanyam2020k}. A company decides to invest in a subset of $n$ projects, each with an uncertain cost $c_i(\bxi)$ and an uncertain profit $r_i(\bxi)$ that both depend on the nominal cost and profit, respectively. Let $\bxi$ be the risk factor that dictates the difference from the nominal values to the actual ones. This risk factor, which we treat as an uncertain scenario, is contained in a given polyhedral uncertainty set $\Xi$. The company has a budget $B$ that it can invest in a project either before or after observing the risk factor $\bxi$. In the latter case, the company generates only a fraction $\eta$ of the profit, which reflects a penalty of postponement. The objective of the capital budgeting problem is to maximize total revenue, resulting in the following formulation:%
\begin{subequations}  \label{eq:CB}%
\begin{align}
    \max_{\vect x \in \gX} \; \min_{\bxi \in \Xi} \; \max_{\vect y \in \gY} \quad & \vect r(\bxi)^\intercal(\vect x + \eta\vect y) \\ 
    \text{s.t.} \quad  & \vect x + \vect y \leq \1 \label{eq:CB_smaller1}\\ 
    & \vect c(\bxi)^\intercal (\vect x + \vect y) \leq B \label{eq:CB_budgetconstraint}.
\end{align}
\end{subequations}
Here, $\gX = \gY = \{0, 1\}^n$ and the $i$-th entries of $\vect x$ and $\vect y$ indicate whether the company invests in the $i$-th project in the first or second stage, respectively. Constraint \eqref{eq:CB_smaller1} ensures that the company invests in a project in at most one of the two stages and constraint \eqref{eq:CB_budgetconstraint} enforces the budget. For each project, the uncertain cost and profit of project $i$ are defined as 
\begin{equation*}
    c_i(\bxi,\bpi) = \big (1 + \bxi\cdot\boldsymbol{\Delta}_i/2 \big ) \Bar{c}_i \quad \text{and} \quad r_i(\bxi,\bpi) = \big (1 +  \bxi\cdot\boldsymbol{\Delta}_{i+n}/2 \big ) \Bar{r}_i, \quad \forall i \in \{1,\ldots,n\},
\end{equation*}
where $\Bar{c}_i$ and $\Bar{r}_i$ are the nominal cost and profit of project $i$, respectively. $\boldsymbol{\Delta}_i$ and $\boldsymbol{\Delta}_{i+n}$ are the $p$-dimensional $i$-th and $(i+n)$-th rows of the parameter matrix  $\boldsymbol{\Delta} \in \R^{2n \times p}$, with $\bxi \in \Xi = [-1, 1]^p$.

\paragraph{Constraint Uncertainty:}
\label{app:cb_trivial}
While uncertain parameters appear in the constraints, i.e., \eqref{eq:CB_budgetconstraint},
ensuring constraint feasibility for \method{} is trivial as for any first-stage decision
$\vect x$ we check if $\max_{\bxi\in \Xi} \vect c(\bxi)^\top \vect x\le B$, where the maximum can be easily calculated since it is a linear problem over $\Xi$. If the latter inequality is true, the second-stage problem is feasible since we can choose $\vect y = \vect 0$.  As such, during the learning augmented CCG, one can check if there is a scenario, $\bxi_0$, that violates the constraint, and if so, add the constraint $\vect c(\bxi_0)^\top \vect x\le B$.  

\paragraph{Instances in Experiments:} The instance parameters are generated similarly as described in \cite{subramanyam2020k}. The uncertainty set dimension $p$ is set to $4$, the nominal cost vector $\boldsymbol{\Bar{c}}$ is chosen uniformly at random from $[0, 10]^n$, $\boldsymbol{\Bar{r}} = \boldsymbol{\Bar{c}}/5$, and $B = \frac{1}{2} \sum_{i=1}^n \Bar{c}_i$. The rows of the sensitivity matrix $\pi$ are sampled uniformly from $[0, 1]^{p}$ and such that each row sums up to 1. This is also known as the unit simplex. We consider instances of sizes $n\in\{10,20,30,40,50\}$.  

\paragraph{Prediction Target:} As discussed by \cite{subramanyam2020k}, the capital budgeting problem has uncertainty in the first- and second-stage objectives, which can easily be transformed into only having second-stage uncertainty. 
For this reason, the prediction target is the sum of the first- and second-stage objectives given the equivalence. 

\subsection{Facility Location Problem}
\label{app:fl}
For facility location, we build on the problem defined in Equation~\ref{eq:FL}.  However, objective uncertainty is also considered in addition to the demand uncertainty.  Specifically, the uncertainty set is defined as $\Xi^o = \{\bxi^o\in[0,1]^n : \sum_{i=1}^n \xi_i^o \leq B^o\}$ with $\xi^o_{i}$ defining a cost uncertainty parameter for facility $i$. This can be understood as a disruption in the facility's operations due to various exogenous factors. The objective coefficients are then given by $\bar d_{ij} + \alpha_i \cdot \xi^o_{i}$, where $\bar d_{ij}$ is the nominal transportation cost from $i$ to $j$, and $\alpha_i$, a deviation factor that regulates the amount by which facility $i$ is impacted by $\xi^o_i$. This formulation is similar in spirit to that of \cite{cheng2021robust_inf}, but with the disruption in operations directly affecting the objective rather than the capacity of the $i$-th facility.  This gives the following problem formulation: 

\begin{subequations}  \label{eq:FL2}%
\begin{align}
    \min_{\vect x \in \gX} \; \max_{\bxi \in \Xi, \bxi^o \in \Xi^o} \; \min_{\vect y \in \gY} \quad & \vect c^\intercal \vect x + \sum_{i=1}^n \sum_{j=1}^m  (\bar d_{ij} + \alpha_{i} \cdot \xi_{i}^o) \cdot y_{ij}  & \\ 
    \text{s.t.} \quad  
        & \sum_{j=1}^m y_{ij} \cdot (\bar{b}_j + \Delta_j \cdot \xi_j)  \leq C_i \cdot x_i, & \forall i \in \{1,\ldots,n\},  \label{eq:FL_budget_constr2}\\ 
        & \sum_{i=1}^n y_{ij} = 1, & \forall j \in \{1,\ldots,m\} \label{eq:FL_demand_constr2}.
\end{align}
\end{subequations}

\paragraph{Instances in Experiments:}
For the facility location problem, we consider instances of size $(n,m) \in \{(5,10), (5,20), (5,50), (10,20), (10,50), (10,20), (10,50),  (20,20), (20,50) \}$, uncertainty budgets $B\in \{ 0.2\cdot m, 0.5\cdot m, 0.8\cdot m\}$ and $B^o \in \{0.1\cdot n, 0.2\cdot n\}$.  For each instance size and uncertainty budget, we randomly sample 5 instances as per the instance generation procedure defined below.

\paragraph{Instance Generation:}
We have two considerations for generating facility location instances.  First, we aim to have instances without recourse, i.e., first-stage decisions will necessarily not have feasible second-stage assignments.  This ensures that the feasibility network can be evaluated.  The second consideration is having problems with discrete second-stage decisions, i.e., a customer can only be assigned to a single facility.  As such, the two-stage robust facility location instances, as specified by \cite{cheng2021robust, cheng2021robust_inf}, do not meet either of these requirements.  Furthermore, even utilizing the base instances to construct problems with binary second-stage decisions leads to relatively trivial instances, as the customers with the largest demands exceed the vast majority of facility capacities, leading to the first-stage decisions that always open the facilities with the largest capacity.  

For these reasons, we utilize the approach outlined by \citet{cornuejols1991comparison} to generate facility location instances, which can easily be adapted to have binary assignment variables.   
Let $I(a,b)$ denote the probability distribution over integers in the range $[a,b]$ that assigns equal mass to each integer in the range.    
Specifically, the location (coordinates) of the customer $j$ and facility $j$ are denoted by ($c_j^x,c_j^y$) and ($f_i^x,f_i^y$), uniformly at random in the unit square.   
The nominal distance from facility $i$ to customer $j$ is given by $d_{ij} = 10 \cdot \sqrt{(f_i^x - c_j^x)^2 + f_i^y - c_j^y)^2}$, i.e., $10$ times the Euclidean distance.  Let $\mathbf{D}\in\mathbb{R}^{n\times m}$ denote the distance matrix with $\mathbf{D}_{ij} = \bar{d}_[ij]$.  Nominal demands (${\bar{b}_j}$) are sampled from $I(5, 35)$, capacities ($C_i$) are sampled from $I(10, 160)$.  The fixed costs ($c_i$) are given by $ I(0, 90) + \sqrt{C_i} \cdot I(100, 110)$.  Capacities are then rescaled as $C_i = C_i \cdot 4 \cdot \frac{\sum_{j=1}^m \bar{b}_j}{\sum_{k=1}^n C_k}$.  To construct the uncertainty sets, we follow an approach similar to that of \cite{cheng2021robust_inf}.  The demand demand deviations are sampled as $\Delta_j \sim U(0.15 \cdot \bar{b}_j, \bar{b}_j)$ and the objective deviations are sampled as $\alpha_i \sim U(0.15 \cdot \textrm{mean}(\mathbf{D}_{i:}), \textrm{mean}(\mathbf{D}_{i:}))$, where $\mathbf{D}_{i:}$ the $i$-th row of $\mathbf{D}$.

\section{Training Data Sampling}
\label{app:sampling}
In this section, we elaborate on the discussion in~\Cref{sec:sampling} and specify how scenarios are sampled for each problem.  Firstly, consider the general approach for budgeted uncertainty sets, i.e., $\Xi = \{\bxi\in[0,1]^q : \sum_{j=1}^q \xi_j \leq B \}$.  In this case, we sample each index of $\bxi$ independently, i.e., $\xi_i \sim U(0,1)$, and normalize to sum to the sampled budget $\bxi = b \cdot \bxi / \sum_{i=1}^q \xi_i$.  After normalization if there exists an index $i$, such that $\xi_i > 1$, then we set $\xi_i = 1$.  The following sections refer to this as sampling for budgeted uncertainty sets (SBUS).

\begin{itemize}
    \item \textbf{Knapsack}:  We sample $b \sim U(\frac{1}{B}, {B})$, then sample using SBUS with budget $b$. 

    \item \textbf{Capital Budgeting}:  We sample each index $\xi_i \sim U(-1,1)$.

    \item \textbf{Facility Location}:  We separately sample the uncertainty for demand and objective uncertainty.   

    \begin{itemize}
    \item  For $\Xi$, we sample in two ways with equal probability.  For the first approach, we use SBUS with a budget of $B$.  For the second approach, we sample binary demands, i.e., assigning all the $\bxi$ to a subset of customers.  Specifically, we define a probability distribution based on the demands and deviations via softmax, i.e., $\xi_j^d$ is set to 1 with probability $\sigma_j = \frac{e^{\bar b_j + \Delta_j}}{ \sum{k=1}^{q_1} e^{\bar b_k + \Delta_k}}$, without replacement, for $\lfloor B \rfloor$ indices of $\bxi$.
    \item  For $\Xi^o$, we can first observe that the worst-case scenarios will always disrupt only the set of open facilities.  For this reason, we sample in two ways with equal probability.
    For the first approach, we sample using SBUS with budget $B^o$.
    For the second approach, we sample using SBUS with budget $B^o$, but only for the same subset of open facilities, i.e., $x_i = 1$.
    \end{itemize}  
\end{itemize}

\section{Objective Value Evaluation}
\label[appendix]{app:ove}
 When we compare the calculated solutions of \method{} and the baseline in our experiments, we need to calculate the objective value of a solution $\vect x^\star\in \mathcal X$ exactly or approximately. The former involves solving the adversarial problem \eqref{eq:adversarial_appendix} for a given solution. Solving this problem is intractable when we have uncertain parameters in the constraints. We first expand on how the adversarial would be solved in a tractable way if the uncertain parameters only appear in the objective function. Subsequently, we describe an approach to approximately solve the AP, which is based on sampling scenarios from $\Xi$.

\subsection{Objective uncertainty}
For the special case of objective uncertainty, the adversarial problem can be solved much more efficiently. In this case, the adversarial problem is given as
\begin{subequations}
\begin{align}
    \max_{\bxi\in\Xi} \min_{\vect y\in\mathcal{Y}} & 
        \quad  \vect c^\intercal \vect x^\star + \vect d(\bxi, \bpi)^\intercal \vect y \\
    \text{s.t.} & 
        \quad W \vect y \le \vect h - T \vect x^\star, &
\end{align}
\label{eq:adversarial_objectiveuncertainty_appendix}%
\end{subequations}
which can be reformulated as
\begin{subequations}
\begin{align}
    \max_{\bxi\in\Xi} \quad & \alpha \\
    \text{s.t.} \quad & \alpha \le \vect c^\intercal \vect x^\star + \vect d(\bxi, \bpi)^\intercal \vect y  \quad \forall \vect y\in \bar{\mathcal Y}, 
\end{align}
\label{eq:adversarial_reformulation_objectiveuncertainty_appendix}%
\end{subequations}
where $\bar{\mathcal Y}=\left\{ \vect y\in \mathcal Y: \quad W \vect y \le \vect h - T \vect x^\star\right\}$. While the set $\bar{\mathcal Y}$ can contain an exponential number of solutions, the latter problem can be solved by iteratively generating the constraints for $y\in \bar{\mathcal Y}$.

\subsection{Constraint uncertainty}
We collect all scenarios $\xi\in \Xi$ generated during the baseline algorithm's and our algorithm's solution procedures.
In addition, we sample 10,000 scenarios of which the details are described in~\Cref{app:sampling}.  
Combining these, we get a set of scenarios $\Xi^{samples}$ and for a first-stage decision $\vect x$ solve the problem 
\begin{subequations}
\begin{align}
    \max_{\bxi\in\Xi^{samples}} \min_{\vect y\in\mathcal{Y}} & 
        \quad \vect c^\intercal \vect x^\star + \vect d(\bxi, \bpi)^\intercal \vect y \\
    \text{s.t.} & 
        \quad W(\bxi,\bpi) \vect y \le \vect h(\bxi,\bpi) - T(\bxi,\bpi) \vect x^\star, & \quad%
\end{align}
\label{eq:adversarial_sample_appendix}%
\end{subequations}
The latter problem can be solved by calculating the optimal value of the second-stage problem for each scenario independently and choosing the worst-case overall optimal values. 

\section{Machine Learning Details}
\label{app:ml_details}
This section reports training details and features. For all problems, the same architecture from Figure~\ref{fig:full_architecture}, with varying parameters, is trained. However, for knapsack and capital budgeting, the feasibility network, $\Phi_6$, is not required as feasibility is either guaranteed or trivial to enforce.  
{In addition, since the facility location problem contains multiple sources of uncertainty, i.e., demand ($\bxi$) and objective ($\bxi^o$), two separate sets of scenario embedding networks are trained.  We denote the networks as $\Phi_3$ and $\Phi_4$ to embed $\bxi$, and $\Phi_3^o$ and $\Phi_4^o$ to embed $\bxi^o$.}
The sizes of all networks are given in~\Cref{app:ml_params}. Generally, all networks have one layer, ranging from 16 to 64 nodes for the embedding networks ($\Phi_1, \ldots, \Phi_4)$ and 8 for the prediction networks ($\Phi_5$ and $\Phi_6$). 

A single dataset of 250,000 (for knapsack and capital budgeting) or 500,000 (for facility location) data points and a model are deployed for all instances of each problem. For knapsack, 36 minutes are spent on data collection and 63 minutes on training. For capital budgeting, 53 minutes on data collection and 37 minutes on training. For facility location, 50 minutes on data generation and 3.7 hours for training. Overall, the longest training ``preprocessing'' time is $\sim 4.5$ hours, comparable to the $3$ hour time limit for solving each instance. Since only a single dataset and model is trained for \textit{all} instances of a benchmark, the amortized cost is relatively insignificant to that of solving. In the following, more details on the features, hyperparameters of the network, and training curves are given, per benchmark.

\subsection{Features}
\label{app:ml_feats}
Table~\ref{tab:ml_feats} reports the features.  For all problems, the set-based architecture, as defined in Figure~\ref{fig:full_architecture} is utilized, so we report the features for each dimension of first-stage decision and uncertainty.  

\begin{table*}[htbp!]\centering\resizebox{0.8\textwidth}{!}{
\begin{tabular}{l|l|l}
\toprule
Problem           & First-Stage Features & Scenario Features   \\
\midrule
Knapsack            & $x_i, f_i, \bar{p}_i, \hat{p}_i, r_i, c_i, t_i, C$ &  $\xi_i, f_i, \bar{p}_i, \hat{p}_i, r_i, c_i, t_i, C$ \\[0.1cm]
Capital budgeting   & $x_i, r_i, c_i$ &  $\big (1 + \boldsymbol{\Phi}_i^\intercal \bxi/2 \big )_i, \big (1 + \boldsymbol{\Psi}_i^\intercal \bxi/2 \big )_i, r_i, c_i$ \\[0.2cm]
\multirow{2}{*}{Facility Location} & \multirow{2}{*}{$x_i, f_i^x, f_i^y, C_i, c_i, \alpha_i, \sum_{i=1}^n C_i, \text{stats}(\mathbf{D}_{i:}) $} & 
        demand ($\bxi$): $ \xi_j^d, c_j^x, c_j^y, \bar{b}_j, \Delta_j, \sum_{i=1}^n C_i, \text{stats}(\mathbf{D}_{:j}), \text{stats}(\alpha) $ \\ 
    & & objective ($\bxi^o$): $ \xi_i^o, f_i^x, f_i^y, C_i, c_i, \alpha_i, \sum_{i=1}^n C_i, \text{stats}(\mathbf{D}_{i:}) $ \\
\bottomrule
\end{tabular}}
\caption{Features for first-stage decision variables and scenarios.  For facility location, the features for each type of uncertainty are separated.  The features $\textrm{stats}(\mathbf{u})$ refer to a tuple of statistics that are computed over the list $\mathbf{u}$, namely, the $\textrm{min}(\mathbf{u})$, $\textrm{max}(\mathbf{u})$, $\textrm{mean}(\mathbf{u})$, and $\textrm{median}(\mathbf{u})$. 
}
\label{tab:ml_feats}
\end{table*}

\subsection{Hyperparameters}
\label{app:ml_params}
Table~\ref{tab:ml_params} reports the hyperparameters for each problem.  As the objective of \method{} is to enable efficient optimization, we train small networks that can achieve a low mean absolute error value to ensure that the main and adversarial problems are tractable. For this reason, no systematic hyperparameter tuning was done. Hyperparameter optimization would likely only further improve the already strong numerical results.  For both problems, we train a model for 500 epochs and compute the mean absolute error on a validation set every 10 epochs. We then use the model with the lowest reported mean absolute validation error during training for evaluation.  As facility location has value and feasibility prediction, we first train the model for 250 epochs, then fix the first-stage and scenario embedding networks, i.e., $\Phi_1,\Phi_2, \Phi_3, \Phi_4, \Phi_3^o, \Phi_4^o$, and only update the parameters for the value and feasibility networks $\Phi_5$ and $\Phi_6$.

\begin{table*}[htbp!]\centering\resizebox{0.5\textwidth}{!}{
\begin{tabular}{l|l|l|l}
\toprule
Hyperparameter                          & Knapsack  & Capital budgeting  & Facility Location         \\
\midrule
Feature scaling                         & min-max   & min-max         & min-max         \\
Label scaling                           & min-max   & min-max          & min-max       \\
\# epochs                               & 500       & 500              & 500     \\
Batch size                              & 256       & 256              & 256       \\
Learning rate                           & 0.001     & 0.001           & 0.001        \\
Dropout                                 & 0         & 0              & 0   \\
Loss function                           & MSELoss   & MSELoss            & MSELoss       \\
Optimizer                               & Adam      & Adam            & Adam         \\
$\Phi_1$  dimensions      & [32, 16]  & [16, 4]            & [16, 4]       \\
$\Phi_2$ dimensions             & [64, 8]   & [32, 8]            & [32, 8]           \\
$\Phi_3$ dimensions          & [32, 16]  & [16, 4]            & [16, 4]       \\
$\Phi_4$ dimensions                & [64, 8]   & [32, 8]           & [32, 8]        \\
$\Phi_3^o$ dimensions          & -  & -            & [16, 4]       \\
$\Phi_4^o$ dimensions                & -   & -           & [32, 8]        \\
$\Phi_5$ dimensions                       & [8]       & [8]               & [8]   \\
$\Phi_6$ dimensions                       & -       & -               & [8]   \\
Aggregation type                        & sum       & sum             & sum    \\
\bottomrule
\end{tabular}}
\caption{Hyperparameters for neural networks.  A value of - for $\Phi_6$ indicates no feasibility prediction is required for the problem in the case of $\Phi_6$
In addition, $\Phi_3^o$ and $\Phi_4^o$ are only required for the facility location problem as there are two types of uncertainty.  
}
\label{tab:ml_params}
\end{table*}

\subsection{Training Curves}
\label{app;ml_training}
\begin{figure}[htbp!]
    \centering
    \subfloat[Knapsack Problem]{
        \includegraphics[width=0.6\textwidth]{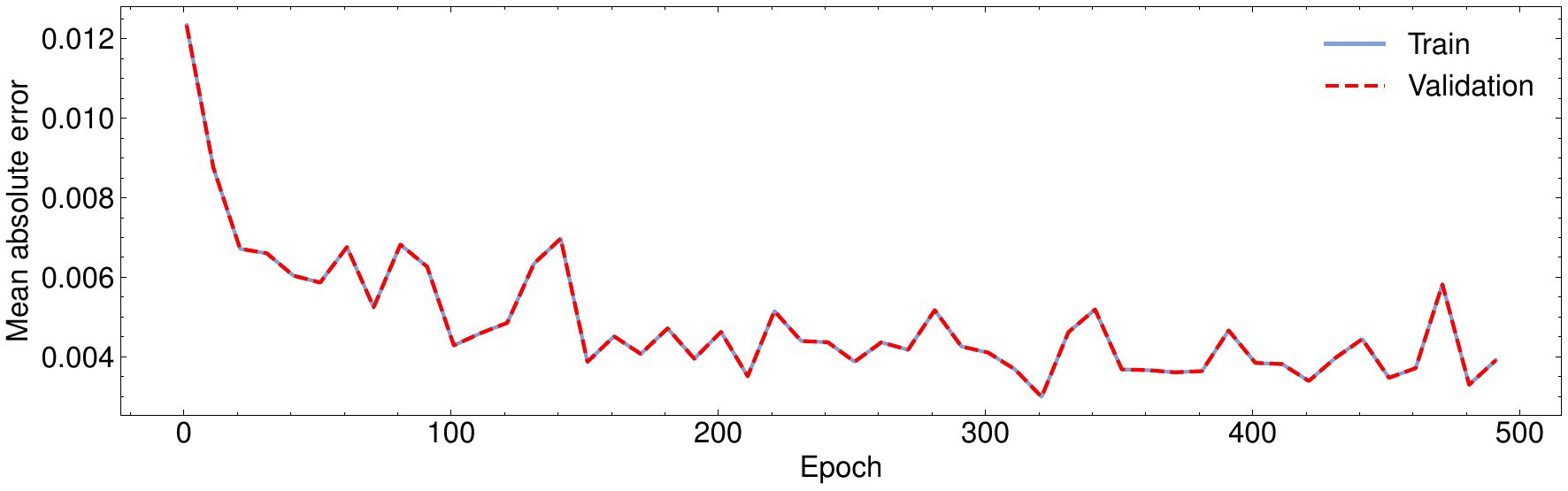}} \\
    \subfloat[Capital Budgeting Problem]{
        \includegraphics[width=0.6\textwidth]{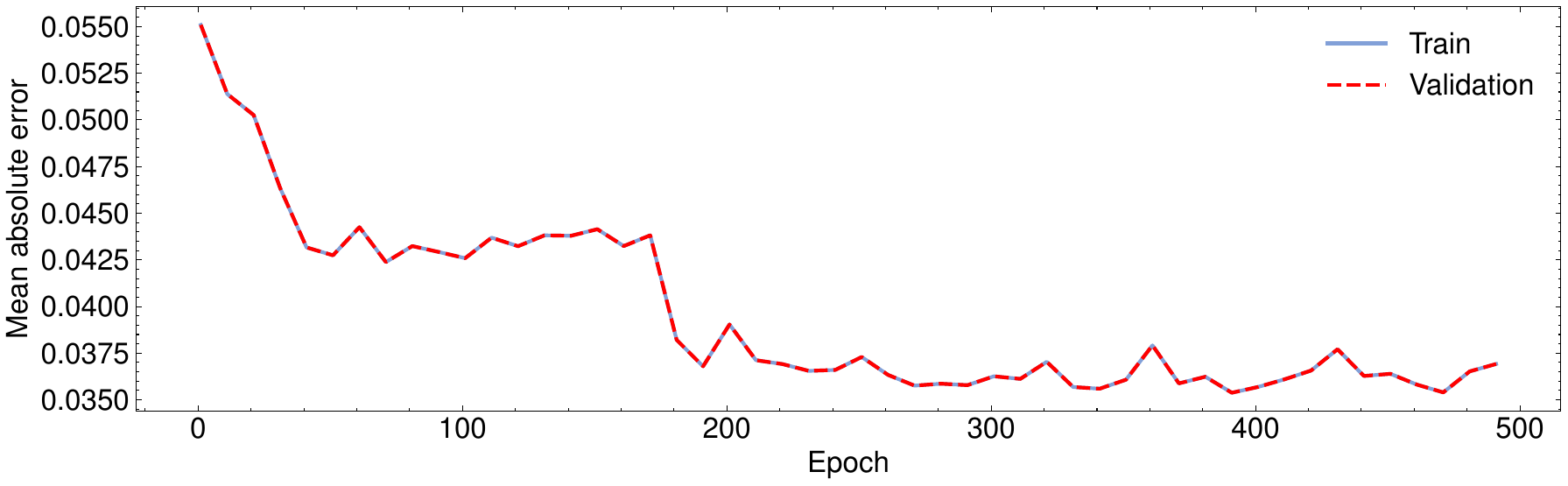}} \\
    \subfloat[Facility Location Problem]{
        \includegraphics[width=0.6\textwidth]{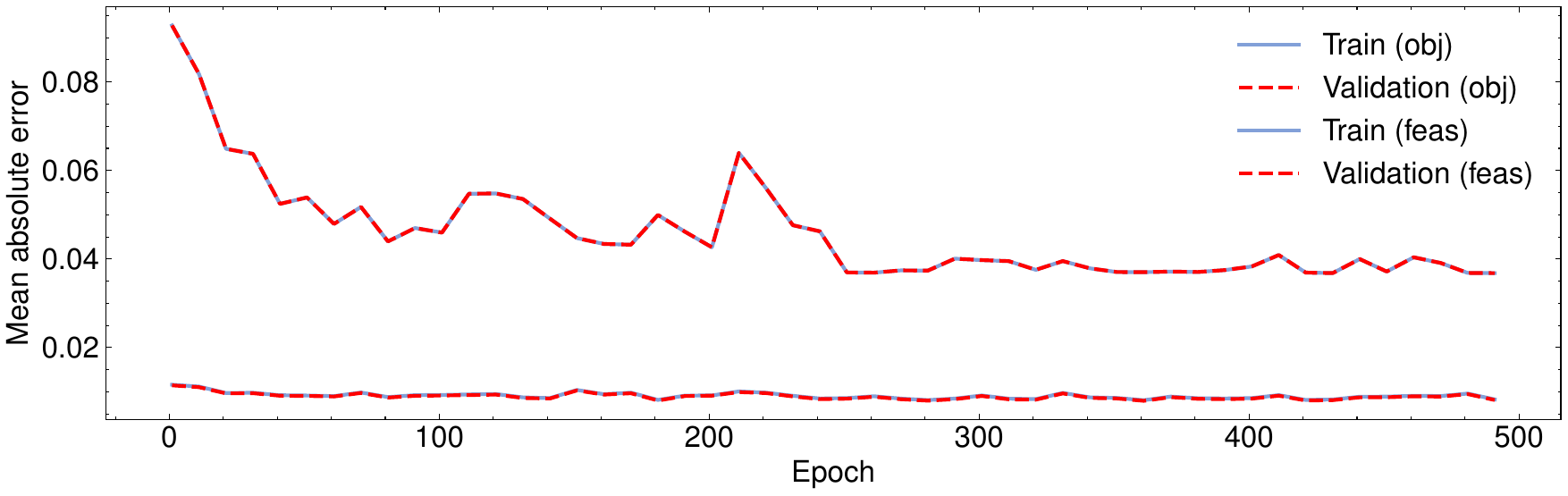}}
    \caption{Training and validation loss over training epochs for Knapsack, Capital Budgeting, and Facility Location.}
    \label{fig:training_curves}
\end{figure}

Figure~\ref{fig:training_curves} reports the validation performance over epochs for each problem.  

\section{Ablation Study}
\subsection{Main Problem Formulation}
\label{app:ab_argmax}
As an alternative to the formulation using $\arg\max$ over a set of scenarios, one can formulate using the $\max$ scenarios. For this ablation study, we only perform this alternative formulation on the knapsack problem which does not have constraint uncertainty. Therefore, the main problem is formulated as:
\begin{subequations}
\begin{align}
    \min_{\vect x\in\mathcal{X}, \alpha} 
        &  \quad  \vect c^\intercal \vect x  + \alpha &  \label{eq:max_obj} \\
    \text{s.t.}  & \quad \alpha \ge   \text{NN}^O_{\Theta}(\bpi, \vect x, \bxi),  \qquad \forall \bxi \in \Xi'_O.  \label{eq:max_cons}
\end{align}
\end{subequations}
Table~\ref{tab:kp_results_maxvsargmax} reports the median relative of the $\arg\max$ and $\max$ formulations and the solving time.  Table~\ref{tab:kp_results_maxvsargmax} demonstrates an improvement in solution quality, with $\arg\max$ obtaining a lower median relative error in every case and a lower computing time in most cases.

\begin{table}[htbp!]\centering\resizebox{1.0\textwidth}{!}{
    \centering
    \small
    \begin{tabular}{l|c|rr|rrcl|c|rr|rr}
        \cmidrule[0.8pt]{0-5}
        \cmidrule[0.8pt]{8-13}
        Correlation  & \# items & \multicolumn{2}{c|}{Median Relative Error}    &  \multicolumn{2}{c}{Times}   & \hspace{0.25cm}  &  Correlation  & \# items & \multicolumn{2}{c|}{Median Relative Error}  &  \multicolumn{2}{c}{Times}  \\
        Type         &          & $\arg\max$  & $\max$                   &  $\arg\max$  &  $\max$                       & \hspace{0.25cm}  &  Type         &          & $\arg\max$  & $\max$                    &  $\arg\max$  &  $\max$     \\               
        \cmidrule[0.6pt]{0-5}
        \cmidrule[0.6pt]{8-13}
        \multirow{7}{*}{Uncorrelated}                                                                        
             & 20 & \textbf{0.000} & 1.167 & \textbf{5} & 11                                                & & \multirow{7}{*}{\shortstack[l]{Almost\\ Strongly\\ Correlated}} & 20 & \textbf{0.000} & 2.042 & \textbf{5} & 12 \\
             & 30 & \textbf{0.000} & 0.945 & \textbf{7} & 14                                             & &      & 30 & \textbf{0.000} & 1.433 & \textbf{6} & 14 \\
             & 40 & \textbf{0.000} & 1.931 & \textbf{9} & 24                                            & &      & 40 & \textbf{0.000} & 1.739 & \textbf{11} & 33 \\
             & 50 & \textbf{0.000} & 1.634 & \textbf{10} & 33                                             & &      & 50 & \textbf{0.000} & 3.161 & \textbf{8} & 20 \\
             & 60 & \textbf{0.000} & 0.452 & \textbf{17} & 29                                             & &      & 60 & \textbf{0.000} & 2.449 & \textbf{15} & 30 \\
             & 70 & \textbf{0.000} & 0.801 & \textbf{19} & 28                                             & &      & 70 & \textbf{0.000} & 2.497 & \textbf{18} & 35 \\
             & 80 & \textbf{0.000} & 2.227 & \textbf{13} & 35                                          & &      & 80 & \textbf{0.000} & 1.824 & \textbf{17} & 30 \\
        \cmidrule[0.6pt]{0-5}
        \cmidrule[0.6pt]{8-13}
        \multirow{7}{*}{\shortstack[l]{Weakly\\ Correlated}}
            & 20 & \textbf{0.000} & 3.515 & \textbf{6} & 13                                         & &  \multirow{7}{*}{\shortstack[l]{Strongly\\ Correlated}}   & 20 & \textbf{0.000} & 1.154 & \textbf{5} & 11 \\
            & 30 & \textbf{0.000} & 2.405 & \textbf{11} & 22                                            & &    & 30 & \textbf{0.000} & 0.967 & \textbf{7} & 15 \\
            & 40 & \textbf{0.000} & 0.502 & \textbf{26} & 42                                          & &    & 40 & \textbf{0.000} & 1.928 & \textbf{16} & 28 \\
            & 50 & \textbf{0.000} & 0.254 & \textbf{24} & 39                                         & &    & 50 & \textbf{0.000} & 3.613 & \textbf{10} & 21 \\
            & 60 & \textbf{0.000} & 1.528 & 77 & \textbf{58}                                         & &    & 60 & \textbf{0.000} & 2.005 & \textbf{20} & 26 \\
            & 70 & \textbf{0.000} & 1.769 & \textbf{18} & 35                                        & &    & 70 & \textbf{0.000} & 2.657 & \textbf{16} & 33 \\
            & 80 & \textbf{0.000} & 3.492 & \textbf{27} & 75                                          & &    & 80 & \textbf{0.000} & 2.051 & \textbf{16} & 28 \\
        \cmidrule[0.8pt]{0-5}
        \cmidrule[0.8pt]{8-13}
    \end{tabular}}
    \caption{Median relative error and solving times for knapsack instances for $\arg\max$ and $\max$ formulations on knapsack instances.  For each row, the median RE and average solving time are computed over 18 instances. All times in seconds. The smallest (best) values in each row/metric are in bold. }
    \label{tab:kp_results_maxvsargmax}
\end{table}

\subsection{Performance of Baselines at \method{} Termination Time}
\label{app:ml_times}
In this section, we report the solution quality, i.e., the median relative error, for Static Robust Optimization (\static{}) and $k$-adaptability at the termination time of \methodmip{} in Tables~\ref{tab:cb_results_ml_time} and \ref{tab:cflp_results_ml_time}.  For capital budgeting, we observe that \methodmip{} still achieves better quality solutions across the same set of instances, and in fact, the relative error is consistently lower than the results at 3 hours, demonstrating the efficacy of \methodmip{}, even if baselines are given a similar runtime.  For facility location, solution quality is similar.  However, baselines such as $k$-adaptability are unable to even compute a feasible solution within the time \methodmip{} takes to terminate.  For the knapsack problem, we compare directly to results reported in \cite{arslan2022decomposition}, so evaluating the solution quality over time is not possible.  

\begin{table*}[htbp!]\centering\resizebox{0.4\textwidth}{!}{
\begin{tabular}{c|rrrrr}
\toprule
\# items & \multicolumn{5}{c}{Median Relative Error}  \\
 & \methodmip{} & \static{} & $k=2$ & $k=5$ & $k=10$  \\
\midrule
10 & 0.816 & 7.198 & 1.194 & 0.180 & \textbf{0.121}  \\
20 & \textbf{0.000} & 0.669 & 0.319 & 0.296 & 0.209  \\
30 & 0.085 & 0.411 & 0.134 & 0.132 & \textbf{0.047}  \\
40 & \textbf{0.001} & 0.200 & 0.098 & 0.132 & 0.019 \\
50 & \textbf{0.003} & 0.105 & 0.120 & 0.058 & 0.068  \\
\bottomrule
\end{tabular}}
\caption{Capital budgeting problem median errors at \methodmip{} termination time (\methodipo{} results excluded).  For each row, the results are computed over 50 instances. All times in seconds. The smallest (best) values in each row/metric are in bold.}
\label{tab:cb_results_ml_time}
\end{table*}

\begin{table*}[htbp!]\centering\resizebox{0.8\textwidth}{!}{
\begin{tabular}{l|rrrrr|rrrrr}
\toprule
\# items & \multicolumn{5}{c|}{Median RE} & \multicolumn{5}{c}{Percent of feasible/found solutions}  \\
 & \methodmip{} & Static & $k=2$ & $k=5$ & $k=10$ & \methodmip{}\ \  & Static\ \  & $k=2$\ \  & $k=5$\ \  & $k=10$\ \  \\
\midrule
(5, 10) & \textbf{0.000} & \textbf{0.000} & \textbf{0.000} & \textbf{0.000} & \textbf{0.000} & 87 & \textbf{100} & \textbf{100} & 97 & 93 \\
(5, 20) & \textbf{0.000} & \textbf{0.000} & \textbf{0.000} & \textbf{0.000} & \textbf{0.000} & 93 & \textbf{100} & \textbf{100} & 93 & 93 \\
(5, 50) & 1.030 & \textbf{0.000} & \textbf{0.000} & \textbf{0.000} & \textbf{0.000} & 87 & \textbf{100} & 47 & 50 & 47 \\
(10, 10) & \textbf{0.000} & 2.877 & 1.162 & \textbf{0.000} & \textbf{0.000} & 93 & \textbf{100} & 47 & 17 & 13 \\
(10, 20) & \textbf{0.000} & 4.530 & 3.830 & 5.757 & \textbf{0.000} & \textbf{100} & \textbf{100} & 13 & 10 & 3 \\
(10, 50) & \textbf{0.000} & 2.057 & - & - & - & 83 & \textbf{100} & - & - & - \\
(20, 20) & \textbf{0.000} & 7.228 & - & - & - & 80 & \textbf{100} & - & - & - \\
(20, 50) & \textbf{0.000} & 6.736 & - & - & - & 93 & \textbf{100} & - & - & - \\
\bottomrule
\end{tabular}}
\caption{Facility location median errors at \methodmip{} termination time (\methodipo{} and \methodngc{} results excluded).  For each row, the results are computed over 30 instances. All times in seconds. The smallest (best) values in each row/metric are in bold.  A value of ``-" indicates no feasible solutions are found across any instance.}
\label{tab:cflp_results_ml_time}
\end{table*}

\section{Static Robust Optimization (SRO) Baseline}
\label{app:static}
This section contains the problem formulations of the static robust optimization baseline for each of the three benchmark problems. The static formulation of a 2RO problem is simply turning the ``wait-and-see'' into ``here-and-now'' variables such that the problem structure does not allow for adaptability to the scenario anymore; hence, it is static.  

\begin{subequations}
\begin{align}
    \min_{\vect x\in\gX, \vect y \in \gY} \max_{\bxi\in\Xi} \quad & 
        \vect c^\intercal \vect x + \vect d(\bxi)^\intercal \vect y \\
    \text{s.t.} \quad & T(\bxi)\vect x + W(\bxi) \vect y \le \vect h(\bxi), 
\end{align}
\end{subequations}
or equivalently,
\begin{subequations}
\begin{align}
    \min_{\vect x\in\gX, \vect y \in \gY} \quad & \theta \\ 
        \text{s.t.} \quad & \vect c^\intercal \vect x + \vect d(\bxi)^\intercal \vect y \leq \theta, & \forall \bxi \in \Xi, \\
    & T(\bxi)\vect x + W(\bxi) \vect y \le \vect h(\bxi), & \forall \bxi \in \Xi,  
\end{align}
\label{eq:static}%
\end{subequations}
A computationally tractable formulation of \eqref{eq:static} can be derived by including the dual form of each uncertain constraint, see \cite{ben2009robust}.

\subsection{Knapsack Problem}
The uncertainty set of this problem is formulated as~$\Xi = \{\bxi \in [0, 1]^n: \sum_{i=1}^n \xi_i \leq \Gamma\}$.
\begin{subequations}
    \begin{align*}
    \min_{\substack{\vect x \in \{0, 1\}^n, \theta \in \R\\ \vect y \in \{0, 1\}^n, \vect r \in \{0, 1\}^n \\ \boldsymbol{\zeta}^1 \in \R^n_{\geq 0}, \zeta^2 \in \R_{\geq 0}}} \quad & \theta \\
    \text{s.t.} \quad     & \sum_{i=1}^n [(f_i - \Bar{p}_i)x_i - f_i y_i + \zeta^1_i] + \zeta^2 \Gamma \leq \theta, \\
    & \Hat{p}(y_i - r_i) - \zeta^1_i - \zeta^2 \leq 0, & \forall i \in \{1, \ldots, n\} \\
    & \sum_{i=1}^n c_i y_i + t_i r_i \leq C, \\ 
    & r_i \le y_i \leq x_i, & \forall i \in \{1, \ldots, n\}.
\end{align*}
\end{subequations}

\subsection{Capital Budgeting Problem}
The uncertainty set is formulated as~$\Xi = [-1, 1]^p$, where $p=4$ in the instances.
\begin{subequations}%
\begin{align*}
    \max_{\substack{\vect x \in \gX, \vect y \in \gY, \theta \in \R \\ \boldsymbol{\zeta}^1 \in \R^p_{\geq 0}, \boldsymbol{\zeta}^2 \in \R^p_{\geq 0} \\ \boldsymbol{\zeta}^3 \in \R^p_{\geq 0}, \boldsymbol{\zeta}^4 \in \R^p_{\geq 0} } } \quad 
    & \theta \\ 
    \text{s.t.} \quad & \sum_{i=1}^n \bar{r}_i(x_i + \eta y_i) - \sum_{j=1}^p (\zeta^1_j - \zeta^2_j )\geq \theta, \\
    & \frac{1}{2} \sum_{i=1}^n \bar{r}_i (x_i + \eta y_i) \Delta_{i+n, j} - \zeta^1_j + \zeta^2_j = 0, & \forall j \in \{1, \ldots, p\}, \\
    & x_i + y_i \leq 1, & \forall i \in \{1, \ldots, n\}, \\ 
    & \sum_{i=1}^n \bar{c}_i(x_i + y_i) + \sum_{j=1}^p (\zeta^3_j - \zeta^4_j )\leq B, \\
    & \frac{1}{2} \sum_{i=1}^n \bar{c}_i (x_i + y_i) \Delta_{i, j} + \zeta^3_j - \zeta^4_j = 0, & \forall j \in \{1, \ldots, p\}. \\
\end{align*}
\end{subequations}

\subsection{Facility Location Problem}
The uncertainty set is given as~$\Xi = \{\bxi \in {[0, 1]}^m : \sum_{j=1}^m \xi_j \leq B\}$
\begin{subequations}%
\begin{align*}
    \min_{\substack{\vect x \in \gX, \vect y \in \gY, \theta \in \R \\
    \delta^1 \in \R_{\geq 0}, \; \boldsymbol{\delta}^2 \in \R^n_{\geq 0} \\
    \boldsymbol{\zeta}^{(i, 1)} \in \R^m_{\geq 0}, \; \forall i \in \cup [n] \\
    \zeta^{(i, 2)} \in \R_{\geq 0}, \; \forall i \in \cup [n] }} 
    \quad & \theta & \\ 
    \text{s.t.} \quad  
        & \vect c^\intercal \vect x + \sum_{i=1}^n \sum_{j=1}^m \bar{d}_{ij} y_{ij} + \sum_{i=1}^n \delta^1_i + B \delta^{2} \leq \theta, \\
        & \alpha_{i} \sum_{j=1}^m y_{ij} - \delta^{1}_i - \delta^{2} \leq 0 & \forall i \in \{1, \ldots, n\}\\
        & \sum_{j=1}^m (y_{ij} \bar{b}_j + \zeta^{(i, 1)}_j) + B \zeta^{(i, 2)} \leq C_i \cdot x_i, & \forall i \in \{1, \ldots, n\},\\
        & \Delta_j y_{ij} - \zeta^{(i, 1)}_j - \zeta^{(i, 2)} \leq 0, & \forall j \in \{1, \ldots, m\}, \forall i \in \{1, \ldots, n\},\\ 
        & \sum_{i=1}^n y_{ij} = 1, & \forall j \in \{1, \ldots, m\}.
\end{align*}
\end{subequations}

\end{document}